\documentclass[11pt,a4paper]{article}
\usepackage[T1]{fontenc}
\usepackage[utf8]{inputenc}
\usepackage{lmodern,microtype}
\usepackage[margin=27mm]{geometry}
\usepackage{amsmath,amssymb,amsthm,mathtools,booktabs}
\usepackage{tikz}
\usepackage{hyperref,bookmark,xurl}
\hypersetup{colorlinks=true,linkcolor=blue,citecolor=blue,urlcolor=blue,
pdftitle={Maximal sections of the cube},pdfauthor={Vu Khac Ky}}
\newtheorem{theorem}{Theorem}[section]
\newtheorem{proposition}[theorem]{Proposition}
\newtheorem{lemma}[theorem]{Lemma}
\newtheorem{corollary}[theorem]{Corollary}
\theoremstyle{remark}\newtheorem{remark}[theorem]{Remark}
\numberwithin{equation}{section}
\begin{document}
\title{Maximal sections of the cube}
\author{Vu Khac Ky\\[3pt]\small FPT University, Hanoi, Vietnam}
\date{23 November 2026}
\maketitle
\begin{abstract}
We prove the folklore conjecture on maximal central sections of the cube
in every dimension and codimension. The maximizing sections are exactly
the products of diagonals of disjoint coordinate blocks whose sizes
differ by at most one. This determines the sharp volume in the cases
left open by Ball's 1989 bounds. In fixed dimensions, the logarithmic
volume deficit controls both the direction and the displacement of an
affine section linearly, with optimal exponent. We also obtain the sharp
uniform bound on projection densities for independent real random
variables with bounded densities.
\end{abstract}
\noindent\textit{2020 Mathematics Subject Classification.}
52A40, 52A38; 26D15, 60D05.\\
\textit{Key words.} Cube sections, Brascamp--Lieb inequality, tight frames,
relative entropy, stability, projection densities.
\tableofcontents
\section{Introduction}
\label{sec:introduction}

Which \(k\)-dimensional section through the center of a cube has the
largest volume? A longest line section joins opposite vertices.
A largest hyperplane section takes a diagonal in two coordinates and
leaves the others free. Both examples partition the coordinates into
blocks and use one diagonal per block. Among sections of this form,
the volume is largest when the block sizes differ by at most one. The
folklore conjecture on maximal cube sections asserts that no other
subspace gives a larger volume~\cite[Conjecture~1]{NayarTkocz2023}.
We prove this conjecture and determine every case of equality. This
resolves the general extremal problem that Ball studied nearly four
decades ago~\cite[p.~252]{Ball1989}.

Let $Q_n=[-1/2,1/2]^n$, and let $G_{n,k}$ be the Grassmannian of
$k$-dimensional linear subspaces of $\mathbb R^n$. We study
\[
 M(n,k)=\max_{H\in G_{n,k}}|Q_n\cap H|_k,\qquad 1\le k\le n,
\]
where $|\cdot|_k$ denotes Euclidean $k$-dimensional volume.
Write $n=qk+r$, with $0\le r<k$. Partition the coordinates into
$k-r$ blocks of size $q$ and $r$ blocks of size $q+1$, and take the
span of their normalized diagonal vectors. The resulting section is an
orthogonal box. Its side lengths are the square roots of the block sizes,
so its volume is
\begin{equation}\label{eq:balanced-target}
 W_{n,k}=q^{(k-r)/2}(q+1)^{r/2}.
\end{equation}
Indeed, the square of the volume of any such box is the product of
its block sizes. If two sizes \(b,c\) satisfy \(c\ge b+2\), transferring
one coordinate from the larger to the smaller increases their product
by \(c-b-1\). This explains the two consecutive sizes in
\eqref{eq:balanced-target}.
We call these the balanced coordinate-block sections, allowing arbitrary
coordinate signs and permutations. Their subspaces form a finite family
$\mathcal B_{n,k}$.

\begin{theorem}[Maximal sections of the cube]\label{thm:full-balanced-value}
For every pair of integers $n\ge k\ge1$,
\[
 M(n,k)=W_{n,k}.
\]
Equality holds precisely for $H\in\mathcal B_{n,k}$.
\end{theorem}

The proof of Theorem~\ref{thm:full-balanced-value}, including attainment
and the characterization of equality, has been verified in Lean~4.\footnote{The
\href{https://github.com/vukhacky/max-section-cube/blob/163d411e85a27017b58d0e7861993d2559b9c6fc/lean/CubeSections/Main.lean}{formal proof}
and the \href{https://github.com/vukhacky/max-section-cube/tree/163d411e85a27017b58d0e7861993d2559b9c6fc/verification}{verification records}
are available at \url{https://github.com/vukhacky/max-section-cube}. The formalization covers the main theorem.
The stability estimates and applications are proved in
Sections~\ref{sec:stability} and~\ref{sec:consequences}.}

For example, in \(Q_7\) the three-dimensional maximizers have side
lengths \(\sqrt2,\sqrt2,\sqrt3\) and volume \(2\sqrt3\). Ball's
continuous upper bound in this case is \((7/3)^{3/2}\).
The distinction comes from the integer sizes of the three blocks.

\subsection{The extremal section problem}

The connection with probability predates the extremal problem.
For \(U\) uniform on the unit cube and a unit vector \(a\), the
density of the weighted sum \(a\cdot U\) at zero is the volume of
the section perpendicular to \(a\). P\'olya's
1913 work on products of sinc functions is an early instance of this
connection~\cite{Polya1913,NayarTkocz2023}.
In 1979 Hensley~\cite{Hensley1979} bounded the volumes of central hyperplane
sections independently of dimension and conjectured the sharp upper
bound $\sqrt2$. Ball proved this conjecture in
1986~\cite{Ball1986}. In the same period, Vaaler~\cite{Vaaler1979}
proved Good's conjecture that every central section of the unit cube
has volume at least one.
The higher-codimension problem was already under study in 1987, when
Sbai Msahli proposed the sharp bound \(2^{(n-k)/2}\) for
\(n\le2k\)~\cite{Sbai1987}. In 1989 Ball treated arbitrary
codimension~\cite{Ball1989}. His two bounds give
\[
 M(n,k)\le
 \min\left\{(n/k)^{k/2},\,2^{(n-k)/2}\right\}.
\]
The first bound is attained when \(k\) divides \(n\), and the second
when \(n\le2k\). For nondivisible pairs with \(n>2k\), both bounds are
strictly larger than the volume of a balanced block section.

The precise balanced formula~\eqref{eq:balanced-target} appears as
Conjecture~3.14 in Pooley's 2012 dissertation~\cite{Pooley2012}.
Ivanov subsequently formulated the geometric conjecture that every
maximizer is an affine cube~\cite{Ivanov2021Cross}. Ivanov and
Tsiutsiurupa derived variational conditions for a maximizing section
and proved the conjecture for every planar section~\cite{IvanovTsiutsiurupa2021}. Nayar and Tkocz record the general
statement as a folklore conjecture in their
survey~\cite[Conjecture~1]{NayarTkocz2023}.
The new range for the sharp value is \(k\ge3\), \(n>2k\), and
\(k\nmid n\). Here the maximizing box must use two distinct integer
block sizes. The discrepancy with Ball's continuous bound can be large.
With
\(t=r/k\), the logarithm of the ratio between the continuous bound and
\(W_{n,k}\) is \(k/2\) times the difference between \(\log(q+t)\) and its
chord on \([0,1]\). Since \((q+1)^{-2}\le-(\log(q+t))''\le q^{-2}\),
the interpolation remainder gives
\[
 \frac{kt(1-t)}{4(q+1)^2}\le
 \log\frac{(n/k)^{k/2}}{W_{n,k}}
 \le\frac{kt(1-t)}{4q^2}.
\]
At fixed \(q\), if \(t\) stays away from zero and one, the two volumes
therefore differ by a factor exponential in the section dimension.

\subsection{Proof and consequences}

The continuous bound treats the total squared length of the coordinate
projections as freely divisible. At a block section each coordinate
belongs to just one of \(k\) directions, so the multiplicities are
integers. We need a volume estimate that detects the angles between
these directions before optimizing their lengths.

An isometry onto a section writes the cube constraints as slabs in
\(\mathbb R^k\). For a block section, each pair of normals is parallel
or orthogonal. The following estimate measures the cost of intermediate
angles. Let \(K\) be an intersection of slabs of widths \(w_i\) and
positive volume, with unit normals \(u_i\) and positive weights satisfying
\(\sum_i a_i u_i u_i^{\mathsf T}=I_k\). If \(\theta_{ij}\) is the angle
between \(u_i\) and \(u_j\), then
\[
 \sum_{i,j}a_i a_j\cos^2\theta_{ij}\sin^2\theta_{ij}
 \le 20\log\frac{\prod_i w_i^{a_i}}{|K|}.
\]
The logarithm on the right is the deficit in the geometric
Brascamp--Lieb inequality. Ball's argument uses the inequality to
bound the volume. We keep a nonnegative determinant term in Barthe's
transport proof and use it to control the angles. The same identity
couples the one-dimensional transports to a single Gaussian vector. The diagonal of the Gram projection specifies
how much of that cost can be assigned to each pair. Theorem~\ref{thm:incidence} gives the resulting pairwise bound.

The second step recovers the integer block sizes. Duplicating or deleting a
row confines its squared length to the interval between two consecutive
integer reciprocals at a maximizer of the volume ratio. Replacing a row
by a copy of another gives a packing constraint on the squared entries
of the Gram projection. For \(n\ge5k\) that constraint and the angular
bound already account for the entire difference between the continuous
volume bound and the balanced value.

For the three remaining quotients with \(n\ge2k\), we express the
entropy excess from the row lengths and the angular deficit using the
same pair weights. The weights of two classes restrict the angle between their normals through a
minor of the complementary projection. A second restriction comes from
replacing one row from each class by a multiple of their sum. The
determinant increases to first order as the new strip narrows, while the replacement cuts off two
triangles whose area is quadratic in the displacement. Fixing the two
coordinates leaves a section of dimension \(k-2\), so induction bounds
the volume lost over those triangles. Section~\ref{sec:extremal-comparison}
derives these two restrictions. Together they reduce the three quotients
to a comparison in the two squared row lengths. Appendix~\ref{sec:comparisons}
gives a finite rational verification, common to all nineteen pairs of
multiplicities.

Two stability statements follow. For cube sections, let \(P_H\) be
the orthogonal projection onto \(H\), and let \(d(H,\mathcal B_{n,k})\)
be the minimum Hilbert--Schmidt distance from \(P_H\) to a balanced
block projection. For \(n>k\) and every positive-volume affine section
with \(y\in H^\perp\), we prove
\[
 d(H,\mathcal B_{n,k})+|y|
 \le C_{n,k}\log\frac{W_{n,k}}{|Q_n\cap(H+y)|_k}.
\]
The volume has a first-order loss at every block section, even an
unbalanced one. This gives the optimal exponent and explicit local
constants. A planar example shows why the global constant must depend
on the dimensions.

The extremal formula also answers a question about concentration.
If independent real random variables have densities bounded by one,
how large can the density of a \(d\)-dimensional projection be?
The linear-image inequality of Madiman, Melbourne, and Xu reduces
this question to uniform variables~\cite{MMX2017}. For the uniform
law on the cube the maximum density is exactly the central section
volume perpendicular to the projection. Consequently the sharp
constant is \(W_{n,n-d}\), and a projection that nearly attains it must
have a kernel close to one of the balanced subspaces
(Corollary~\ref{cor:product-marginals}).

Two other consequences use different parts of the proof.
The angular estimate controls the cost of moving the normal directions
to the axes of one orthonormal basis. Applied to supporting slabs,
it gives a quantitative volume inequality for convex bodies
(Corollary~\ref{cor:widths}). Polarity gives a further application of
the sharp section formula. Combining it with the three-dimensional
symmetric Mahler inequality determines the minimum volume of a
three-dimensional projection of the cross-polytope
(Corollary~\ref{cor:polar-minimum}).

\section{Frames, determinant entropy, and Gaussian transport}
\label{sec:frames}

We first write the section as an intersection of slabs in Euclidean
space. Isotropy makes its volume accessible to the geometric
Brascamp--Lieb inequality. The frame formulation and its variational
use are established in~\cite{Ivanov2021Frames,IvanovTsiutsiurupa2021}. Our purpose is to
retain information about the angles that is lost in the usual scalar
bound.

Choose an isometric parameterization $A\in\mathbb R^{n\times k}$
of $H$, so $A^{\mathsf T}A=I_k$. For each nonzero row write
$A_i=\sqrt{c_i}u_i^{\mathsf T}$, with $|u_i|=1$. Then
\[
 \sum_i c_i u_i u_i^{\mathsf T}=I_k,
 \qquad \sum_i c_i=k,
 \qquad 0<c_i\le1.
\]
The numbers $c_i$ are the squared lengths of the coordinate projections,
often called the row leverages. Thus $c_i=|P_H\varepsilon_i|^2$ for the standard coordinate vectors
$\varepsilon_i$ of $\mathbb R^n$. The section is isometric to
\[
 K=\{x\in\mathbb R^k:|u_i\cdot x|\le r_i\text{ for all nonzero rows}\},
 \qquad r_i=\frac1{2\sqrt{c_i}}.
\]
Zero rows impose no constraint and are omitted from frame sums. They
will be included with $c_i=0$ in scalar entropy sums.

The angular defect is
\[
 \mathcal E=\sum_{i,j}c_i c_j(u_i\cdot u_j)^2
                          \{1-(u_i\cdot u_j)^2\}.
\]
It is nonnegative and vanishes precisely when every pair of nonzero
row directions is parallel or orthogonal. In that case isotropy forces
the distinct unoriented directions to form an orthonormal basis.

The Jacobian in the transport argument is a weighted sum of
rank-one matrices. Cauchy--Binet expresses its determinant as an
average over bases. We record the corresponding probability law,
which will also give a way to select nearly orthogonal directions
in Section~\ref{sec:axes}. Let $a_i>0$, $|u_i|=1$, and
$\sum_i a_i u_i u_i^{\mathsf T}=I_k$. Put $P_i=u_i u_i^{\mathsf T}$
and let $U$ have columns $\sqrt{a_i}u_i$. For every $k$-element index
set $S$, define
\[
 p(S)=\det(U_S)^2.
\]
Cauchy--Binet and $UU^{\mathsf T}=I_k$ imply $\sum_Sp(S)=1$.
Only linearly independent sets have positive probability. The generating
identity
\begin{equation}\label{eq:core-dpp-generating}
 \sum_{|S|=k}p(S)\prod_{i\in S}t_i
 =\det\left(\sum_i a_i t_iP_i\right)
\end{equation}
gives the inclusion probabilities
\begin{equation}\label{eq:core-dpp-inclusions}
 \mathbb P(i\in S)=a_i,
 \qquad
 \mathbb P(i,j\in S)=a_i a_j\{1-(u_i\cdot u_j)^2\}\quad(i\ne j).
\end{equation}
Indeed, the mixed derivative of $\det(I+sX+tY)$ at zero is
$\operatorname{tr}X\operatorname{tr}Y-\operatorname{tr}(XY)$.

For a finite potential vector $z=(z_i)$, set
\begin{equation}\label{eq:core-determinant-remainder}
 J(z)=\log\det\left(\sum_i a_i e^{z_i}P_i\right)-\sum_i a_i z_i.
\end{equation}
If \(S\) is sampled with probabilities \(p(S)\), write
\(Y=\sum_{i\in S}z_i\). Cauchy--Binet gives
\(J(z)=\log\mathbb E e^Y-\mathbb E Y\), which is nonnegative by Jensen's
inequality. Equivalently, reweighting each basis by \(e^Y/\mathbb E e^Y\)
gives a probability law \(p_z\) with \(J(z)=D_{\rm KL}(p\|p_z)\).
Here \(D_{\rm KL}(\nu\|\lambda)=\int\log(d\nu/d\lambda)\,d\nu\)
is relative entropy when \(\nu\ll\lambda\). The law \(p\) is a
projection determinantal process~\cite{Lyons2003}.

\subsection{An exact Gaussian transport identity}

The geometric Brascamp--Lieb inequality belongs to a general theory
of multilinear integral inequalities, with Gaussian extremization and
structural decompositions developed in~\cite{BrascampLieb1976,Lieb1990,BCCT2008}.
We use the increasing transports from Barthe's proof~\cite{Barthe1998}
and retain the terms in its change of variables. This gives an exact
decomposition of the deficit. Relative entropy enters because it
compares the transported law with the standard Gaussian. The duality
between entropy and Brascamp--Lieb inequalities is developed more
generally by Carlen and Cordero-Erausquin~\cite{CarlenCordero2009}.

For the moment take any isotropic unit frame $(a_i,u_i)$ and positive
half-widths $r_i$. Let $K$ be their open slab intersection, let $V=|K|$,
and put
\[
 B=\prod_i(2r_i)^{a_i}.
\]
Write $\Phi$, $\varphi$, and $\gamma_k$ for the standard Gaussian
distribution function, one-dimensional density, and $k$-dimensional law.
The coordinates of a uniform point of the intersection need not be
uniform on their full intervals. We nevertheless use the maps that
take each full interval with uniform measure to the standard Gaussian.
Their derivatives then encode the slab widths exactly. Define
\[
 T_i(s)=\Phi^{-1}\left(\frac{s+r_i}{2r_i}\right),
 \qquad z_i(x)=\log T_i'(u_i\cdot x),
 \qquad F(x)=\sum_i a_iT_i(u_i\cdot x)u_i.
\]
Thus $\varphi(T_i(s))T_i'(s)=1/(2r_i)$.

Let $X$ be uniform on $K$. Set
\[
 b(X)=(\sqrt{a_i}T_i(u_i\cdot X))_i,
 \quad Q=I-U^{\mathsf T}U,
 \quad P(X)=|Qb(X)|^2,
 \quad \mu=F_\#\operatorname{Unif}(K).
\]
Since $UU^{\mathsf T}=I_k$, $Q$ is an orthogonal projection and
$F(X)=Ub(X)$. In particular
$P(X)=\sum_i a_iT_i(u_i\cdot X)^2-|F(X)|^2$.

\begin{lemma}\label{lem:core-entropy-coupling}
With $\delta_{\rm BL}=\log(B/V)$, one has
\begin{equation}\label{eq:core-entropy-decomposition}
 \delta_{\rm BL}
 =D_{\rm KL}(\mu\|\gamma_k)+\mathbb E J(z(X))
                                     +\tfrac12\mathbb E P(X).
\end{equation}
There is a coupling of $X$ and $G\sim\gamma_k$ such that
\begin{equation}\label{eq:core-lifted-coupling}
 \mathbb E\sum_i a_i
       |T_i(u_i\cdot X)-u_i\cdot G|^2
 \le2\{\delta_{\rm BL}-\mathbb E J(z(X))\}
 \le2\delta_{\rm BL}.
\end{equation}
\end{lemma}
\begin{proof}
The differential
$DF(x)=\sum_i a_iT_i'(u_i\cdot x)P_i$ is positive definite.
The map is the gradient of a strictly convex function on the convex
set $K$, and is therefore injective. The change-of-variables formula
on its image gives
\[
 \int_K\varphi_k(F(x))\det DF(x)\,dx\le1.
\]
Using \eqref{eq:core-determinant-remainder} and the one-dimensional
transport identities, its integrand is exactly
\begin{equation}\label{eq:core-transport-integrand}
 \frac1B\exp\{J(z(x))+P(x)/2\}.
\end{equation}
Both terms in the exponent are nonnegative.

The density change under the injective map also gives
$d\mu/dy(F(x))=1/(V\det DF(x))$. Therefore
\[
 \log\frac{d\mu}{d\gamma_k}(F(x))
 =\delta_{\rm BL}-J(z(x))-P(x)/2.
\]
All terms are integrable. The body $K$ is bounded, its one-dimensional
marginals have bounded densities by a finite fiber-volume bound, and
$\int_{-r_i}^{r_i}T_i(s)^2\,ds=2r_i$. Thus the transports and $F(X)$
have finite second moments. Integrability of $J$ and $P$ follows from
\eqref{eq:core-transport-integrand}. Integrating the logarithmic identity
proves \eqref{eq:core-entropy-decomposition}.

Talagrand's Gaussian transportation inequality \cite{Talagrand1996},
in the normalization $\operatorname{Cov}(\gamma_k)=I_k$, states
\[
 W_2(\mu,\gamma_k)^2\le2D_{\rm KL}(\mu\|\gamma_k).
\]
Take a quadratic-cost coupling of $F(X)$ and $G$ with this bound,
and lift it to $X$. The two vectors in
\[
 b(X)-U^{\mathsf T}G
   =Qb(X)+U^{\mathsf T}(F(X)-G)
\]
are orthogonal pointwise. Consequently their squared norms add, and
\[
 \mathbb E|b(X)-U^{\mathsf T}G|^2
 \le\mathbb E P+2D_{\rm KL}(\mu\|\gamma_k)
 =2\{\delta_{\rm BL}-\mathbb E J\}.
\]
This is \eqref{eq:core-lifted-coupling}.
\end{proof}

The decomposition in particular proves the geometric Brascamp--Lieb
bound $V\le B$. For cube rows, $B=\prod_{c_i>0}c_i^{-c_i/2}$.
Concavity of entropy and $\sum_i c_i=k$ then give
\[
 V\le B\le(n/k)^{k/2}.
\]
This is the continuous bound used at the endpoints of the final
maximum-ratio argument.

\section{The determinant remainder and the angles between normals}
\label{sec:exponential-determinant-coercivity}

The transport identity bounds the total deficit from below by the
expected determinant remainder. We need to know which directions are
responsible for it. If the potentials of two rows differ, their
interaction should disappear when the rows are orthogonal. The
following inequality makes this precise, with the squared inner product
as the weight of each pair. Its second part gives a bounded function
of the potential that we can compare under the Gaussian coupling.

For an isotropic unit frame $(a_i,u_i)$, write
\[
 J(z)=\log\det\left(\sum_i a_i e^{z_i}u_i u_i^{\mathsf T}\right)
                      -\sum_i a_i z_i,
 \qquad
 \mathcal V(v)=\sum_{i<j}a_i a_j\rho_{ij}^2(v_i-v_j)^2,
 \quad \rho_{ij}=u_i\cdot u_j.
\]
Isotropy gives the identity
\[
 \mathcal V(v)=\sum_i a_i v_i^2-
 \left\|\sum_i a_i v_i u_i u_i^{\mathsf T}\right\|_{\rm HS}^2.
\]
Hence $0\le\mathcal V(v)\le\sum_i a_i v_i^2$.

\begin{lemma}[A pairwise lower bound for the determinant remainder]
\label{lem:pairwise-determinant-production}
For every real potential vector $z$,
\begin{equation}\label{eq:pairwise-determinant-production}
 J(z)\ge\sum_{i<j}a_i a_j\rho_{ij}^2
                 \phi(|z_i-z_j|),\qquad \phi(t)=t-1+e^{-t}.
\end{equation}
If every $z_i\ge0$, then
\begin{equation}\label{eq:exponential-determinant-production}
 J(z)\ge\frac92\mathcal V((e^{-z_i/3})_i).
\end{equation}
\end{lemma}
\begin{proof}
A common additive shift leaves both sides of
\eqref{eq:pairwise-determinant-production} unchanged, so assume
$\min_i z_i=0$. Raise a common cutoff from zero. While the cutoff
moves, only the coordinates above it change, so a pair contributes to
the derivative exactly when the cutoff lies between its two potentials.
Put $z_i(t)=\min(z_i,t)$ and, away from the finitely
many switching times, let
\[
 A_t=\{i:z_i>t\},\qquad
 S_t=\sum_{i\in A_t}a_i u_i u_i^{\mathsf T},\qquad
 X_t=e^{-t}\sum_i a_i e^{z_i(t)}u_i u_i^{\mathsf T}.
\]
The matrix $X_t$ is positive definite and satisfies $X_t\preceq I_k$.
Differentiating the determinant gives
\[
 \frac{d}{dt}J(z(t))=\operatorname{tr}\{S_t(X_t^{-1}-I_k)\}.
\]
For the pairwise expression on the right side of
\eqref{eq:pairwise-determinant-production}, only pairs with one index
in $A_t$ contribute to the derivative.  Since $\phi'(s)=1-e^{-s}$,
that derivative is
\[
 \sum_{i\in A_t,j\notin A_t}a_i a_j\rho_{ij}^2(1-e^{z_j-t})
 =\operatorname{tr}\{S_t(I_k-X_t)\}.
\]
Their difference is nonnegative, because
\[
 X_t^{-1}-2I_k+X_t=(I_k-X_t)X_t^{-1}(I_k-X_t)\succeq0
\]
and $S_t\succeq0$.  Both expressions vanish at $t=0$.
The finitely many switching times do not affect the integral, so integration up to $\max_i z_i$ proves the first assertion.

To transfer this estimate to the Gaussian coupling, we want a bounded
function of the potentials whose squared differences are controlled by
\(\phi\). Exponentiation turns additive changes of potential into
changes of scale. The first two nonzero terms at zero determine the scale.
\[
 \phi(t)=\tfrac12t^2-\tfrac16t^3+O(t^4),\qquad
 C(1-e^{-\alpha t})^2=C\alpha^2t^2-C\alpha^3t^3+O(t^4).
\]
Matching the two coefficients gives \(\alpha=1/3\) and \(C=9/2\).
These local coefficients give a global bound, because the difference
vanishes at zero and its derivative is \((1-e^{-t/3})^3\ge0\).

Moreover, if $s,t\ge0$, then
\[
 |e^{-s/3}-e^{-t/3}|\le1-e^{-|s-t|/3}.
\]
Apply these two inequalities pair by pair in
\eqref{eq:pairwise-determinant-production} to obtain
\eqref{eq:exponential-determinant-production}.
\end{proof}

\subsection{Two angular bounds from a common Gaussian coupling}
\label{sec:edge-allocation}

We need two consequences of the determinant estimate. The first controls
the total angular defect of a frame. The second distinguishes a row of
small weight from one that carries almost an entire direction. Isotropy specifies how much of
the angular interaction is incident to each row. Write
\(s_{ij}=\langle u_i,u_j\rangle^2\) and \(w_{ij}=a_i a_j s_{ij}\).
Taking the quadratic form of the frame identity in direction \(u_i\)
gives
\begin{equation}\label{eq:pair-degree}
 \sum_{j\ne i}w_{ij}=a_i(1-a_i).
\end{equation}
This identity distributes the cost of transporting the \(i\)th
coordinate among its pairs with the other coordinates.

\begin{theorem}[A weighted angular bound for slab intersections]
\label{thm:incidence}
Let \(\sum_i a_i u_i u_i^{\mathsf T}=I_k\), with unit vectors \(u_i\)
and positive weights. For arbitrary \(r_i>0\), put
\(K=\{x:|\langle u_i,x\rangle|<r_i\ \forall i\}\) and
\(\delta=\log[\prod_i(2r_i)^{a_i}/|K|]\). Define
\[
 M_{\rm G}=\sqrt{3/5},\qquad
 \mathcal L(A,B,s)=
 \frac{M_{\rm G}(1-s)}
 {(16-s)\{2/9+[2/(3e)](2-A-B)\}}.
\]
Then
\begin{equation}\label{eq:incidence-bound}
 \delta\ge\sum_{i<j} w_{ij}\,\mathcal L(a_i,a_j,s_{ij}).
\end{equation}
For all choices of the slab widths, one also has
\[
 \mathcal E:=\sum_{i,j}a_i a_js_{ij}(1-s_{ij})
 \le \frac{32}{\sqrt{3/5}}\left(\frac29+\frac2{3e}\right)\delta
 \le 20\delta .
\]
Both bounds remain valid for intervals with arbitrary centers, with
\(\delta\) computed from their intersection of positive volume.
\end{theorem}
\begin{proof}
Use the common coupling of Lemma~\ref{lem:core-entropy-coupling}.
Write \(T_i=T_i(\langle u_i,X\rangle)\), \(G_i=\langle u_i,G\rangle\),
and \(e_i=T_i-G_i\). Up to a common constant the determinant potentials
are \(z_i=\beta_i+T_i^2/2\), where \(\beta_i=-\log r_i\).
By Lemma~\ref{lem:pairwise-determinant-production} and the coupling,
\[
 \delta\ge\mathbb E\sum_{i<j}w_{ij}\phi(|z_i-z_j|)
                  +\tfrac12\mathbb E\sum_i a_i e_i^2.
\]
Indices with \(a_i=1\) have no pair of positive weight by
\eqref{eq:pair-degree}. Discard their nonnegative transport costs.
For the remaining indices the right side is exactly the sum, over
pairs of positive weight, of \(w_{ij}\) times
\begin{equation}\label{eq:pair-energy}
 \mathbb E\left[\phi(|z_i-z_j|)
       +\frac{e_i^2}{2(1-a_i)}+\frac{e_j^2}{2(1-a_j)}\right].
\end{equation}

Fix such a pair. Subtract \(\min(\beta_i,\beta_j)\) from both potentials,
and set \(t_h=\exp[-(\beta_h-\min(\beta_i,\beta_j))/3]\),
\(h=i,j\). Both \(t_h\) are at most one and one equals one. For
\(f_h(y)=t_h e^{-y^2/6}\) the scalar inequality in the proof of
Lemma~\ref{lem:pairwise-determinant-production} gives
\(\lvert f_i(T_i)-f_j(T_j)\rvert\le\sqrt{2\phi(|z_i-z_j|)/9}\).
Also \(\|f_h'\|_\infty\le t_h/\sqrt{3e}\).
The triangle inequality and Cauchy--Schwarz therefore bound
\(\mathbb E(f_i(G_i)-f_j(G_j))^2\) by \eqref{eq:pair-energy} times
\[
 \frac29+\frac2{3e}\big[(1-a_i)t_i^2+(1-a_j)t_j^2\big]
 \le\frac29+\frac2{3e}(2-a_i-a_j).
\]
For correlated standard Gaussians, integration of
\(e^{-(G_i^2+G_j^2)/6}\) against their joint law gives
\(3/\sqrt{16-s_{ij}}\), while \(\mathbb E e^{-G_i^2/3}=\sqrt{3/5}\).
Thus, putting \(M=\sqrt{3/5}\) and
\(D=3/\sqrt{16-s_{ij}}\), the Gaussian squared difference is
\(M(t_i^2+t_j^2)-2Dt_it_j\). If \(t_i=1\), completion of the square
bounds this expression below by
\[
 M(t_j-D/M)^2+M-D^2/M\ \ge\
 M\,\frac{1-s_{ij}}{16-s_{ij}} .
\]
Interchanging the indices covers \(t_j=1\). Summing proves
\eqref{eq:incidence-bound}. 
For the global bound we sum the squared differences before estimating
the transport errors. The ordinary frame contraction
\(\mathcal V(v)\le\sum_i a_i v_i^2\) also controls the
pair-dependent amplitudes. Write \(t_{i|j}\) for the amplitude of
direction \(i\) in the pair \(i,j\). We claim
\begin{equation}\label{eq:schur-angular}
 \sum_{i<j} w_{ij}(t_{i|j}v_i-t_{j|i}v_j)^2\le\sum_i a_i v_i^2 .
\end{equation}
To see this, put \(b_i=r_i^{-2/3}\) and
\(v_i(x)=v_i\mathbf 1_{[0,b_i]}(x)/\sqrt{b_i}\).
Then \(t_{i|j}=\sqrt{\min(b_i,b_j)/b_i}\). In particular,
\[
 (t_{i|j}v_i-t_{j|i}v_j)^2
 =\int_0^{\min(b_i,b_j)}(v_i(x)-v_j(x))^2\,dx .
\]
Extend every integral to \([0,\infty)\), sum over pairs, and apply the
frame contraction at each \(x\). The resulting upper bound is
\(\sum_i a_i\int_0^\infty v_i(x)^2\,dx=\sum_i a_i v_i^2\).

Put \(j_0=\mathbb E J(z)\), \(h(y)=e^{-y^2/6}\), and let
\(\|D(Y)\|^2=\sum_{i<j}w_{ij}\mathbb E
 [t_{i|j}h(Y_i)-t_{j|i}h(Y_j)]^2\).
The scalar potential estimate gives \(\|D(T)\|^2\le 2j_0/9\).
Apply \eqref{eq:schur-angular} pointwise to \(v_i=h(T_i)-h(G_i)\).
The coupling and \(\|h'\|_\infty^2\le 1/(3e)\) show that the squared
norm of \(D(T)-D(G)\) is at most \(2(\delta-j_0)/(3e)\).
The triangle inequality and Cauchy--Schwarz now give
\[
 \|D(G)\|^2\le
 \left(\sqrt{2j_0/9}+\sqrt{2(\delta-j_0)/(3e)}\right)^2
 \le \left(\frac29+\frac2{3e}\right)\delta .
\]
The Gaussian square completion bounds the same norm below by
\(M_{\rm G}\mathcal E/32\), since \(16-s_{ij}\le16\).
Finally \(M_{\rm G}>77/100\) and \(e>8/3\) give the stated rational
constant. For a translated slab intersection \(K'\), its symmetric
difference body \((K'-K')/2\) lies in the centered intersection with
the same widths. By Brunn--Minkowski, \(|K'|\le |K|\), so the
translated deficit is at least \(\delta\).
\end{proof}

\begin{remark}
One can retain a slightly better coefficient for nearly equal widths.
Use the same amplitude \(d_i=(r_i/\max_jr_j)^{1/3}\) on every
pair incident to \(i\). The determinant and Lipschitz estimates in
the preceding argument still give a total Gaussian cost at most
\([2/9+2/(3e)]\delta\) by the frame contraction
\(\mathcal V(v)\le\sum_i a_i v_i^2\) to
\(v_i=d_i[h(T_i)-h(G_i)]\). The exact Gaussian cost of one pair is
\[
 M_{\rm G}(d_i-d_j)^2+
 2d_i d_j\left(M_{\rm G}-3/\sqrt{16-s_{ij}}\right).
\]
Convexity of \(s\mapsto3/\sqrt{16-s}\) bounds the parenthesis below by
\((M_{\rm G}-3/4)(1-s_{ij})\). Consequently one can also bound
\(\mathcal E\) by
\[
 \frac{2/9+2/(3e)}{M_{\rm G}-3/4}
 \left(\frac{\max_i r_i}{\min_i r_i}\right)^{2/3}\delta .
\]
\end{remark}
\begin{corollary}[Zero deficit and divisible endpoints]
\label{cor:zero-deficit-rigidity}
An intersection of translated slabs of positive volume has zero deficit precisely when its distinct unoriented directions are orthogonal and the slabs in each parallel class coincide. In particular, if an isometric cube-section frame satisfies $V=B_{\rm BL}$, then
its distinct unoriented nonzero row directions are orthogonal.
If $n=qk$ and $V=q^{k/2}$, its parallel classes all have size $q$,
and the section is a balanced coordinate-block section.
\end{corollary}
\begin{proof}
Theorem~\ref{thm:incidence}, with zero deficit,
gives $\mathcal E=0$. Every summand of $\mathcal E$ is nonnegative,
so each pair of nonzero directions is parallel or orthogonal.
The weights in each parallel class sum to one. The length of the
intersection of its intervals is at most the shortest interval length,
hence at most their weighted geometric mean. Equality requires all
lengths to agree and the intervals to coincide, which is also sufficient.
Taking the product over the orthogonal classes proves the slab statement.
At a divisible endpoint, equality in the two inequalities
$V\le B_{\rm BL}\le(n/k)^{k/2}$ is necessary. Strict concavity in
the entropy step, including all $n$ raw rows, gives $c_i=1/q$ for
every row. Isotropy then forces exactly $q$ rows in each orthogonal
class. Coordinate signs account for the possible orientations.
\end{proof}

\section{Integer packing and the reciprocal window}
\label{sec:integer-packing}
\label{sec:integer-leverage-quantization}

The loss in the continuous Brascamp--Lieb bound can be compared directly
with the integer block sizes.  For this purpose a maximizing frame carries
additional information because replacing one row by another cannot
increase the normalized volume.  The resulting two-row
inequality makes every squared Gram row subject to a sharp packing
constraint.

Let $A\in\mathbb R^{n\times k}$ satisfy $A^{\mathsf T}A=I_k$ and write
$P=AA^{\mathsf T}$ and $c_i=P_{ii}>0$.  Thus $P$ is a projection,
$\sum_jP_{ij}^2=c_i$, and $\sum_i c_i=k$.  The angular defect is
\[
 \mathcal E=k-\sum_{i,j}\frac{P_{ij}^4}{c_i c_j}.
\]
We say that the frame satisfies the row-replacement inequalities if
\begin{equation}\label{eq:integer-row-replacement}
 |c_i-c_j|\le c_i c_j-P_{ij}^2\qquad(i\ne j).
\end{equation}
Every volume-maximizing cube-section frame with $c_i<1$ satisfies
these inequalities.  Indeed, replacing row $i$ by a copy of row $j$
removes one constraint from the slab intersection and leaves the other
constraints unchanged.  Its frame operator is
$M=I_k-a_i a_i^{\mathsf T}+a_j a_j^{\mathsf T}$, where $a_i$ denotes
the $i$th row as a column vector.  Normalization by $M^{-1/2}$ produces
another $n$-row isometry with volume at least $V\sqrt{\det M}$.
The operator is positive definite because $c_i<1$.  Maximality gives
$\det M\le1$, and the rank-two determinant formula gives
\[
 \det M=1+c_j-c_i-c_i c_j+P_{ij}^2.
\]
Interchanging $i,j$ proves \eqref{eq:integer-row-replacement}.

\begin{lemma}[Packing a squared Gram row]
\label{lem:integer-gram-packing}
Suppose \eqref{eq:integer-row-replacement} holds.  Put
$q_i=\lfloor1/c_i\rfloor$.  Then
\[
 \boxed{\quad
 \mathcal E\ge
 \sum_i(1-q_i c_i)\bigl((q_i+1)c_i-1\bigr).
 \quad}
\]
The summands are nonnegative, and vanish when $c_i^{-1}$ is an integer.
In a common reciprocal window $c_i\in[1/(q+1),1/q]$, the right side can
be written
\[
 T_q(c)=\sum_i(1-qc_i)\bigl((q+1)c_i-1\bigr).
\]
At the lower endpoint the expression with $q_i=q+1$ and the expression
with $q$ both vanish.
\end{lemma}
\begin{proof}
If $c_i\ge c_j$, the replacement inequality implies
\[
 P_{ij}^2\le c_i c_j-c_i+c_j
 =c_j^2+(c_i-c_j)(c_j-1)\le c_j^2.
\]
Together with the same statement with the indices exchanged, and the
diagonal identity, this gives
\[
 0\le P_{ij}^2\le\min(c_i^2,c_j^2).
\]
The numbers $Y_{ij}=P_{ij}^2$ in row $i$ have sum $c_i$ and are at most
$c_i^2$.  Among numbers in $[0,c_i^2]$ with this sum, the sum of their
squares is maximized by filling $q_i$ entries to $c_i^2$ and placing
the remainder in one further entry.  This follows by moving mass
from a smaller positive, nonfull entry to a larger nonfull entry.
Consequently
\[
 \sum_jY_{ij}^2\le q_i c_i^4+(c_i-q_i c_i^2)^2.
\]
Symmetry of $Y$ and the arithmetic--geometric mean inequality give
\[
 \sum_{i,j}\frac{Y_{ij}^2}{c_i c_j}
 \le \frac12\sum_{i,j}Y_{ij}^2(c_i^{-2}+c_j^{-2})
 =\sum_i c_i^{-2}\sum_jY_{ij}^2.
\]
Subtracting from $\sum_i c_i=k$ yields
\[
 \mathcal E\ge\sum_i\{c_i-q_i c_i^2-(1-q_i c_i)^2\}
 =\sum_i(1-q_i c_i)((q_i+1)c_i-1).
\]
This proves the lemma.
\end{proof}

\begin{remark}[An exact nonnegative remainder]
\label{rem:integer-exact-remainder}
In a common reciprocal window define
\[
 \mathcal P_q=
 \sum_i\left\{q c_i^2+(1-qc_i)^2
                     -\frac1{c_i^2}\sum_jP_{ij}^4\right\},
 \qquad
 \mathcal H=
 \frac12\sum_{i,j}P_{ij}^4(c_i^{-1}-c_j^{-1})^2.
\]
Both are nonnegative by the preceding proof, and the following identity
holds exactly:
\begin{equation}\label{eq:integer-exact-production-remainder}
 \boxed{\mathcal E=T_q(c)+\mathcal P_q+\mathcal H.}
\end{equation}
The first remainder measures the failure to fill each squared Gram row
with as many full entries as possible.  The second records differing
leverages along nonzero Gram entries.  These terms are retained jointly;
no choice of an orthonormal basis or assignment to clusters is needed.
\end{remark}

The packing estimate has produced a quadratic that vanishes when
the reciprocal squared row length is an integer. To convert this
estimate into a volume bound, consider the entropy function
\(-c\log c\) on the same interval. At its two endpoints the
corresponding contributions to the balanced volume are exact.
Subtracting the chord through these endpoints leaves a gap controlled
by the same quadratic.

\begin{lemma}[The exact integer entropy gap]
\label{lem:integer-entropy-gap}
Let $n=qk+r$, $0\le r<k$, and suppose $c_i\in[1/(q+1),1/q]$
and $\sum_i c_i=k$.  Set
\[
 B=\prod_i c_i^{-c_i/2},\qquad
 W=q^{(k-r)/2}(q+1)^{r/2},\qquad
 A_q=(q+1)\log(1+1/q)-1.
\]
Then
\begin{equation}\label{eq:integer-entropy-quantization}
 0\le\log(B/W)\le\frac{A_q}{2}T_q(c).
\end{equation}
Here the scalar entropy excess is
\(\psi_q(c)=-c\log(qc)-(1-qc)\log(1+1/q)\).

The coefficient $A_q$ is optimal in the pointwise bound for $\psi_q$.  If at least one leverage is strictly inside
the reciprocal window, the upper inequality is strict.
\end{lemma}
\begin{proof}
Write $a=1/(q+1)$ and $b=1/q$.  The chord of the concave function
$f(c)=-c\log c$ between $a,b$ is
\[
 L_q(c)=c\log q+(1-qc)\log(1+1/q).
\]
Since $\sum_i c_i=k$ and $n-qk=r$,
\[
 \sum_iL_q(c_i)=k\log q+r\log(1+1/q)=2\log W.
\]
Thus \(2\log(B/W)=\sum_i\psi_q(c_i)\), where \(\psi_q=f-L_q\ge0\).
It remains to bound the gap above the chord. Set
\(T(c)=(c-a)(b-c)/(b-a)=(1-qc)((q+1)c-1)\). The quadratic
\(A_q T(c)\) vanishes at both endpoints and has the same derivative
as \(\psi_q\) at \(a\), namely \(A_q\). The quadratic interpolation
remainder therefore gives, at every \(a<c<b\),
\[
 \psi_q(c)-A_qT(c)=\frac{\psi_q^{(3)}(\xi)}{6}(c-a)^2(c-b)<0
 \qquad (a<\xi<b),
\]
since \(\psi_q^{(3)}(\xi)=1/\xi^2\). The ratio \(\psi_q/T\) tends
to \(A_q\) at \(a\), proving sharpness. Summing over the rows gives
the assertion.

\end{proof}

\begin{corollary}[Volume comparison from the angular deficit]
\label{cor:integer-direct-criterion}
Assume that a cube-section frame satisfies the reciprocal window and
row-replacement inequalities, and put
$\delta_{\rm BL}=\log(B/V)$.  If
$\mathcal E\le C\delta_{\rm BL}$, then
\[
 \log(W/V)\ge
 \left(1-\frac{CA_q}{2}\right)\delta_{\rm BL}
 +\frac{A_q}{2}(\mathcal P_q+\mathcal H).
\]
In particular, $CA_q/2<1$ implies $V\le W$, with strict inequality
whenever $\delta_{\rm BL}>0$.
\end{corollary}
\begin{proof}
Use \eqref{eq:integer-entropy-quantization} and
\eqref{eq:integer-exact-production-remainder} in
$\log(W/V)=\delta_{\rm BL}-\log(B/W)$.
\end{proof}

\subsection{The reciprocal leverage window}

Fix \(q\ge2\) and \(k\). Over the \(k+1\) ambient dimensions
\(n=qk,\ldots,(q+1)k\), maximize the ratio \(M(n,k)/W_{n,k}\).
The endpoint ratios are one by the continuous bound and the block
construction. Why vary \(n\) when the goal is an inequality at fixed
\(n\)? Duplicating a row leaves the original slab intersection
unchanged and gives a competitor in dimension \(n+1\). Deleting it
gives one in dimension \(n-1\). At a largest ratio in the interval,
both comparisons are available for each individual row.

\begin{lemma}\label{lem:window}
At such an interior maximum every squared row length satisfies
\(a\le c_i\le b\), where \(a=1/(q+1)\) and \(b=1/q\).
\end{lemma}
\begin{proof}
Let \(V=M(n,k)\). Duplicating row \(i\) preserves the slab body and
multiplies the determinant of the frame operator by \(1+c_i\).
After normalization it gives an \((n+1)\)-row section of volume
\(V\sqrt{1+c_i}\). Since the ratio is maximal and
\(W_{n+1,k}/W_{n,k}=\sqrt{(q+1)/q}\), we get \(c_i\le b<1\).
Deleting the row now leaves full rank and enlarges the slab body.
The normalized volume is at least \(V\sqrt{1-c_i}\), while
\(W_{n-1,k}/W_{n,k}=\sqrt{q/(q+1)}\). This gives \(c_i\ge a\).
\end{proof}

For \(q\ge5\), no information from smaller dimensions is needed.
By Theorem~\ref{thm:incidence} we may take \(C=20\) in
Corollary~\ref{cor:integer-direct-criterion}. Since
\[
 A_q=\int_0^1\frac{1-t}{q+t}\,dt<\frac1{2q},
\]
every maximizer of the volume ratio with \(q\ge5\) satisfies
\begin{equation}\label{eq:large-q-result}
 \log(W_{n,k}/V)\ge (1-10A_q)\,\delta_{\rm BL}\ge0 .
\end{equation}
The first coefficient is strictly positive even at \(q=5\). Equality
therefore forces zero angular defect by Theorem~\ref{thm:incidence}.

\section{The volume comparison at a maximizing frame}
\label{sec:extremal-comparison}

At an interior maximum of the volume ratio, group the nonzero rows by their unoriented
directions. The two replacement inequalities
\eqref{eq:integer-row-replacement} force parallel rows to have the same
length. In reciprocal coordinates they read
\( |c_i^{-1}-c_j^{-1}|\le1-\langle u_i,u_j\rangle^2\), so two directions
with a small angle must have nearly equal reciprocal squared row lengths. Write \(m_i\) for the number of rows in class \(i\), \(c_i\) for
their common squared length and \(a_i=m_ic_i\). The grouped data satisfy
\(\sum_i a_i u_i u_i^{\mathsf T}=I_k\). All rows in class \(i\) impose the
same slab of width \(c_i^{-1/2}\). Keeping it once leaves the intersection
unchanged. Its contribution to the Brascamp--Lieb product is also unchanged,
since \((c_i^{-1/2})^{a_i}=((c_i^{-1/2})^{c_i})^{m_i}\).
Thus Theorem~\ref{thm:incidence} applies with the class weights \(a_i\)
and the original deficit \(\log(B_{\rm BL}/V)\).

A class of weight one is
orthogonal to every other class by \eqref{eq:pair-degree}. Its size
\(1/c_i\) is \(q\) or \(q+1\), and its leverage contributes zero
to the entropy gap. Every class incident to a nonorthogonal pair of distinct classes has
\(a_i<1\) and \(m_i\le q\).

For such classes define
\[
 H_{q,m}(c)=\frac{\psi_q(c)}{2c(1-mc)},\qquad
 \psi_q(c)=-c\log(qc)-(1-qc)\log(1+1/q).
\]
When \(m=q\), the value at \(c=1/q\) is the continuous extension
\((1-q\log(1+1/q))/2\). The degree identity gives the exact comparison
\begin{equation}\label{eq:charges}
 \log(B_{\rm BL}/W_{n,k})
 =\tfrac12\sum_i m_i\psi_q(c_i)
 =\sum_{\substack{i<j\\w_{ij}>0}}w_{ij}
       [H_{q,m_i}(c_i)+H_{q,m_j}(c_j)] .
\end{equation}
Thus the determinant deficit and the entropy gain of nonintegral
leverages are sums with the same weights. A strict inequality on each
pair of positive weight will force the whole section to have volume below the
balanced target.

\subsection{Two coordinates and their fibers}

Fix \(k\ge2\) and \(q\ge2\). Choose an interior ambient dimension
\(qk<n<(q+1)k\) at which the volume ratio \(M(n,k)/W_{n,k}\) is largest
among the interior dimensions in this window. The block construction
gives a ratio at least one, and both endpoint ratios are one. Hence
Lemma~\ref{lem:window} applies. Choose a maximizing section at this
dimension, of volume \(V=M(n,k)\ge W_{n,k}\). Such a section exists by
compactness of the Grassmannian and continuity of the volume, which
follows from the radial function of the cube.

We need two restrictions on the angle between two classes of rows.
One comes from the Gram matrix. The other uses the distribution of the
two coordinates under the uniform law on the section.

For two distinct nonorthogonal classes, let \(m,l\) be their multiplicities and \(c,d\) the
squared lengths of individual rows in the two classes. Orient the rows so the two unit normals
have inner product \(\rho=\sqrt{s}\ge0\).
The grouped frame has Gram projection \(U^{\mathsf T}U\), with
\(U\) having columns \(\sqrt{a_i}u_i\). Positivity of the
corresponding principal minor of \(I-U^{\mathsf T}U\) gives
\begin{equation}\label{eq:complementary-angle}
 s\le\frac{(1-mc)(1-ld)}{mlcd}.
\end{equation}
We next bound the joint density of the two coordinates. Write
\(R=R_q(m,l)=q^{-1}(1+1/q)^{q-(m+l)/2}\).
For the following estimate, when \(k>2\), assume the balanced upper
bound in dimension \(k-2\) for all ambient dimensions at least \(2(k-2)\).
If \(X\) is uniform on the slab body, put
\(Z_i=2\sqrt c\,u_i\cdot X\) and \(Z_j=2\sqrt d\,u_j\cdot X\).
Then
\begin{equation}\label{eq:signed-density-new}
 \text{joint density of }(Z_i,Z_j)\ \le
   \frac{R}{4\sqrt{cd}\sqrt{1-s}}\quad\text{on }[-1,1]^2 .
\end{equation}
Indeed, let \(F\) be the common kernel of the two row functionals.
On \(F\), the \(m+l\) rows in the two selected classes vanish. The remaining
rows define an isometry on \(F\), since the original frame is an
isometry. Translating a fiber to \(F\) therefore identifies it
isometrically with an affine section of \(Q_N\) of dimension \(k-2\),
where \(N=n-m-l\). By Brunn--Minkowski its volume is at most the
central section volume in the same direction. For \(k>2\), the upper
bound in dimension \(k-2\) gives a fiber volume at most \(W_{N,k-2}\).
To see why the fiber stays in the induction range, write \(n=qk+r\).
Since each of the two classes has at most \(q\) rows,
\[
 N-q(k-2)=r+(q-m)+(q-l)\ge r.
\]
Thus removing the two classes preserves the lower bound \(q\) on the
ratio of ambient to section dimension.
When \(k=3\), this bound is simply the diameter bound
\(W_{N,1}=\sqrt N\) for a chord of \(Q_N\).

We compare this bound with \(V\). For every positive integer \(b\),
\[
 b\le q(1+1/q)^{b-q}.
\]
Indeed, the successive ratios \((j+1)/j\) decrease, and equality
holds at \(b=q,q+1\). Multiply the inequality over the \(k-2\) block
sizes that define \(W_{N,k-2}^2\) and divide by
\(W_{n,k}^2=q^k(1+1/q)^{n-qk}\). This gives
\(W_{N,k-2}/V\le W_{N,k-2}/W_{n,k}\le R\).
For \(k=2\) a nonempty fiber is a point. Here \(n=2q+1\), so its
zero-dimensional volume divided by \(V\) is at most
\(1/\sqrt{q(q+1)}\le1/q\le R\).

The coarea Jacobian \(4\sqrt{cd}\sqrt{1-s}\) proves
\eqref{eq:signed-density-new} in either case.

The projection bound can leave the angle unrestricted. For instance,
when \(q=3\), \((m,l)=(1,2)\), and \(c=d=2/7\), its right side
is \(15/8\). We now use maximality to see why two distinct row
directions cannot be arbitrarily close.

\begin{lemma}[Two-row replacement]\label{lem:two-row-comparison}
At a maximizing section, take two distinct nonorthogonal classes
satisfying the density bound~\eqref{eq:signed-density-new}.
Let \(u=\sqrt{1-\rho^2}\), \(h=(c+d)/2\), \(g=\sqrt{cd}\), and
\(L=h+g\rho-g^2u^2\). Then for every \(0\le t\le1/2\),
\begin{equation}\label{eq:exact-strip-law}
 \big[1-h+g\rho+L\{(1-t)^{-2}-1\}\big]
 \left(1-\frac{Rt^2}{gu}\right)_+^2\le1 .
\end{equation}
In particular, writing \(\tau=(1-\rho)/u=\tan(\tfrac12\arccos\rho)\),
\begin{equation}\label{eq:direct-strip-law}
 \frac{L^2}{2R\{1-Lgu/(2R)\}^2}\le \tau+\frac{cd}{8}u^3 .
\end{equation}
\end{lemma}
\begin{proof}
Choose one row \(v\) and one row \(w\) from the two classes, oriented
so that \(v\cdot w=g\rho\). Replace both by \((v+w)/(2(1-t))\),
retaining all other rows. The new constraint keeps precisely those
points of the old section with \(|Z_i+Z_j|\le2(1-t)\).
Figure~\ref{fig:two-corners} shows the complement in the coordinate
square. Its area is \(4t^2\). Multiplication by the density bound
\(R/(4gu)\) shows that the lost fraction of volume is at most
\(Rt^2/(gu)\).

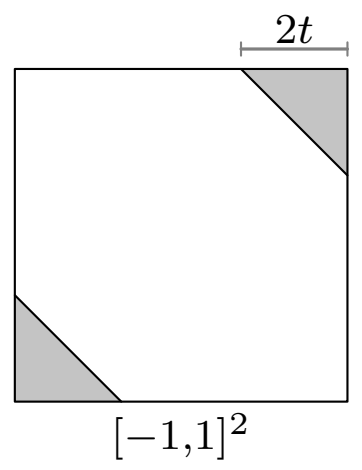
\begin{figure}[htbp]
\centering
\begin{minipage}[c]{.28\linewidth}
\centering
\begin{tikzpicture}[scale=1.10, line cap=round, line join=round]
\fill[black!23] (-1,-1)--(-.36,-1)--(-1,-.36)--cycle;
\fill[black!23] (1,1)--(.36,1)--(1,.36)--cycle;
\draw[thin] (-1,-1) rectangle (1,1);
\draw[thin] (-1,-.36)--(-.36,-1) (.36,1)--(1,.36);
\draw[gray,semithick] (.36,1.12)--(1,1.12);
\draw[gray,thin] (.36,1.08)--(.36,1.16) (1,1.08)--(1,1.16);
\node[above,inner sep=1pt] at (.68,1.12) {$\scriptstyle 2t$};
\node[below,inner sep=2pt] at (0,-1) {$\scriptstyle [-1,1]^2$};
\end{tikzpicture}
\end{minipage}\hfill
\begin{minipage}[c]{.70\linewidth}
\caption{The two-row replacement removes the two shaded triangles
from the coordinate square. Both legs of each triangle have length
\(2t\).}\label{fig:two-corners}
\end{minipage}
\end{figure}

We next normalize the replacement frame. Before normalization its
frame operator is
\[
 I-p_-p_-^{\mathsf T}+\{(1-t)^{-2}-1\}p_+p_+^{\mathsf T},
 \qquad p_\pm=(v\pm w)/\sqrt2 .
\]
It is positive definite since \(|p_-|^2=h-g\rho\le1/2\).
The squared lengths of \(p_\pm\) are \(h\pm g\rho\), and their Gram
determinant is \(g^2u^2\). The rank-two determinant formula therefore
gives the first factor of~\eqref{eq:exact-strip-law}. Normalization
multiplies volume by its square root. The new section cannot have larger
volume than the original one, which proves~\eqref{eq:exact-strip-law}.

At \(t=0\), merging the rows costs \(h-g\rho\) in the determinant.
To choose how far to shrink the strip, write \(\kappa=R/(gu)\).
For \(\kappa t^2<1\), inequality~\eqref{eq:exact-strip-law} becomes
\[
 h-g\rho\ge L\{(1-t)^{-2}-1\}-\{(1-\kappa t^2)^{-2}-1\}.
\]
The leading terms of the determinant gain and the area penalty are
\(2Lt\) and \(2\kappa t^2\), respectively. Their difference is
largest at \(t=L/(2\kappa)=Lgu/(2R)\), which we now use in the
exact bound.
Here \(L\ge h-g^2>0\), \(L\le 2h\le1\), and \(g\le1/q\le R\).
Thus \(t\le1/2\)
and \(Rt^2/(gu)=Lt/2\). Moreover
\[
 (1-Lt/2)^{-2}-1\le\frac{Lt}{1-Lt}\le\frac{Lt}{1-t}
\]
because \((1-Lt/2)^2\ge1-Lt\). Rearranging~\eqref{eq:exact-strip-law}
now gives
\[
 h-g\rho\ge L\{(1-t)^{-2}-1\}-\frac{Lt}{1-t}
          =\frac{Lt}{(1-t)^2}.
\]
Finally, row replacement~\eqref{eq:integer-row-replacement} gives
\(|c-d|\le g^2u^2\). Thus the part of the merger loss due to unequal
row lengths satisfies
\[
 h-g=\frac{(c-d)^2}{4(h+g)}\le\frac{g^3u^4}{8}.
\]
Use \(h-g\rho=g(1-\rho)+(h-g)\) and divide by \(gu\) to obtain
\eqref{eq:direct-strip-law}.
\end{proof}

The scale of the separation is worth recording. If \(\rho\ge1/2\),
then \(L=h+g(\rho-gu^2)\ge h\), since \(g\le1/2\). Also
\(u=(1+\rho)\tau\le2\tau\) and \(cd\le1/4\).
Thus~\eqref{eq:direct-strip-law} gives
\(h^2/(2R)\le\tau+cd\,u^3/8\le 3\tau/2\).
When \(\rho<1/2\), the same conclusion follows from
\(\tau>1/\sqrt3\) and \(h^2/R\le1/q\le1/2\). In either case,
\[
 \tan(\theta/2)\ge\frac{h^2}{3R}.
\]
For fixed \(q\) this keeps the two classes apart independently of the
section dimension. For the quantitative comparison we retain the
stronger bound~\eqref{eq:direct-strip-law}.

\begin{proposition}\label{prop:edge-pay}
Let \(q\in\{2,3,4\}\) and \(1\le m,l\le q\).
Suppose \(c,d\in[1/(q+1),1/q]\) and \(s=\cos^2\theta\) for
\(0<\theta<\pi/2\). If these parameters satisfy
\eqref{eq:complementary-angle} and~\eqref{eq:direct-strip-law},
with \(g,h,L,R\) as above, then
\[
 \mathcal L(mc,ld,s)>H_{q,m}(c)+H_{q,l}(d)+\frac1{4000}.
\]
\end{proposition}
At the chosen maximizer, every pair of distinct nonorthogonal classes
satisfies the hypotheses. Appendix~\ref{sec:comparisons} proves the
proposition by a finite comparison of rational bounds.

\begin{theorem}\label{thm:balanced-at-least-two}
For every \(n\ge2k\), \(M(n,k)=W_{n,k}\), with equality precisely for
the balanced coordinate-block sections.
\end{theorem}
\begin{proof}
We use strong induction on \(k\), proving the upper bound for all
\(n\ge2k\) together. The only lower-dimensional input is the upper
bound for the fibers in~\eqref{eq:signed-density-new}. A line section has length at most \(\sqrt n\), the
diameter of the cube, with equality exactly for its space diagonals.
Now let \(k\ge2\) and fix a quotient \(q\ge2\). At an interior
maximizer in this window, \(V\ge W_{n,k}\). If \(q\ge5\),
\eqref{eq:large-q-result} gives the reverse inequality and forces
\(\delta_{\rm BL}=0\). Hence all distinct classes are orthogonal
by Corollary~\ref{cor:zero-deficit-rigidity}.

For \(q=2,3,4\) the induction hypothesis gives
\eqref{eq:signed-density-new} when \(k>2\), and the point-fiber argument
above gives it when \(k=2\). If two distinct classes were nonorthogonal,
Theorem~\ref{thm:incidence}, Proposition~\ref{prop:edge-pay}, and
\eqref{eq:charges} would give
\[
  \log(B_{\rm BL}/V)>\log(B_{\rm BL}/W_{n,k}),
\]
contrary to \(V\ge W_{n,k}\). Thus in every case the classes are
orthogonal. Each has total weight one by isotropy. Its rows have equal
length by row replacement, so the reciprocal window forces its size
to be \(q\) or \(q+1\). This identifies the section as a balanced block
section, of volume \(W_{n,k}\).

The largest ratio in the window is therefore one. The balanced
construction attains this ratio for each \(n\) in the window, proving
the upper bound throughout. Now take any section of volume \(W_{n,k}\)
at an interior dimension. Its ratio is one, so it too maximizes the
ratio over the whole window. The same argument identifies its signed
coordinate blocks. At the divisible
endpoints the classification follows from
Corollary~\ref{cor:zero-deficit-rigidity}.

\end{proof}
\subsection{Equality in the classical codimension range}
\label{sec:q1-equality}

Write $Q_n=[-1/2,1/2]^n$.  The following argument gives the equality
classification in the classical range $k\le n\le2k$.  Its only additional
analytic inputs are the classical geometric Brascamp--Lieb inequality
for general nonnegative integrable one-dimensional functions, and the
strict form of Ball's scalar integral inequality:
\begin{equation}\label{eq:q1-strict-sinc}
 \int_{\mathbb R}\left|\frac{\sin(\pi t)}{\pi t}\right|^p\,dt
 \le\sqrt{\frac2p}\qquad(p\ge2),
 \quad\text{with equality if and only if }p=2.
\end{equation}
Ball proved the scalar inequality~\cite{Ball1986}. Its strict form
is stated explicitly in~\cite{AbiKhuzam2017}. We use the
standard geometric Brascamp--Lieb inequality~\cite{Barthe1998},
but no classification of its multidimensional equality cases.

\begin{proposition}\label{prop:q1-equality}
Let $k\le n\le2k$, and put $d=n-k$.  A $k$-dimensional linear subspace
$H\subset\mathbb R^n$ satisfies
\[
 \operatorname{vol}_k(Q_n\cap H)=2^{d/2}
\]
if and only if, after a signed permutation of the ambient coordinates,
$H$ is spanned by $d$ mutually disjoint two-coordinate diagonal vectors
and $k-d$ distinct coordinate vectors.  Equivalently, its row classes
are $d$ doubletons and $k-d$ singletons.
\end{proposition}
\begin{proof}
The case \(d=0\) is immediate. Suppose \(1\le d\le k\) and that the
section attains \(2^{d/2}\). We first show that the squared projections
of the coordinate vectors onto \(H^\perp\) satisfy
\[
 b_i:=|P_{H^\perp}\varepsilon_i|^2\le\tfrac12 .
\]
Delete coordinate \(i\) by orthogonally projecting onto the coordinate
hyperplane \(\varepsilon_i^\perp\). On \(H\) the squared volume
Jacobian of this projection is
\(\det(I_H-(P_H\varepsilon_i)(P_H\varepsilon_i)^{\mathsf T})
=1-|P_H\varepsilon_i|^2=b_i\). When \(b_i>0\) the image section
lies in a \(k\)-dimensional section of \(Q_{n-1}\), so Ball's
codimension bound gives
\[
 \sqrt{b_i}\,2^{d/2}\le M(n-1,k)\le 2^{(d-1)/2}.
\]
This proves the claim, including the endpoint \(n-1=k\).

Choose an isometry \(U:\mathbb R^d\to H^\perp\) with rows \(w_i\).
Then \(|w_i|^2=b_i\), \(\sum_i b_i=d\), and the corresponding primal
leverages are \(c_i=1-b_i\). We will use equality in Ball's scalar
bound to show that every positive \(b_i\) equals \(1/2\).

Omit indices with $b_i=0$, and set $v_i=w_i/\sqrt{b_i}$ for the others.
Then $\sum_i b_i v_i v_i^{\mathsf T}=I_d$.  Let
$\operatorname{sinc}(s)=\sin(s)/s$, with its continuous value at zero.
The geometric Brascamp--Lieb inequality yields
\begin{align*}
 \int_{\mathbb R^d}\prod_{i:b_i>0}
       |\operatorname{sinc}(\pi\sqrt{b_i}\,v_i\cdot t)|\,dt
 &\le\prod_{i:b_i>0}
 \left(\int_{\mathbb R}
       |\operatorname{sinc}(\pi\sqrt{b_i}\,s)|^{1/b_i}\,ds
       \right)^{b_i}\\
 &\le\prod_{i:b_i>0}(\sqrt2)^{b_i}=2^{d/2}.
\end{align*}
The last step is~\eqref{eq:q1-strict-sinc} with $p=1/b_i$ and a
one-dimensional change of variable.  In particular, the characteristic
function of $U^{\mathsf T}X$, where $X$ is uniform on $Q_n$, is integrable.
Fourier inversion at zero therefore gives
\[
 \operatorname{vol}_k(Q_n\cap H)
 =\int_{\mathbb R^d}\prod_i
       \operatorname{sinc}(\pi w_i\cdot t)\,dt.
\]
Here the section-volume density furnished by the coarea formula agrees
at zero with the continuous Fourier density because its $1/k$ power is concave
on the interior of the projected cube by Brunn's theorem, so it is
continuous at zero.  Thus the preceding integral bound applies to the
actual central section, not merely to an almost-everywhere density.

If some $0<b_i<1/2$, the corresponding scalar inequality is strict,
and its positive exponent $b_i$ makes the product bound strict.  This
contradicts the assumed equality.  Hence every $b_i$ is either $0$ or
$1/2$, and every $c_i$ is either $1$ or $1/2$.  Since $\sum b_i=d$,
exactly $2d$ of the $c_i$ equal $1/2$.

For the original frame the slab Brascamp--Lieb product is now
\[
 B=\prod_i c_i^{-c_i/2}
   =2^{(2d)/4}=2^{d/2}
   =\operatorname{vol}_k(Q_n\cap H).
\]
Corollary~\ref{cor:zero-deficit-rigidity} implies
that any two row directions are parallel or orthogonal.  Isotropy
makes the sum of the squared raw row lengths in each parallel class
equal to one.  Since each squared length is either $1$ or $1/2$, each
class is either a singleton of squared length $1$ or a doubleton with
both squared lengths $1/2$.  The class directions
form an orthonormal basis, which proves the asserted form.  Conversely,
those block subspaces have $d$ side lengths $\sqrt2$ and $k-d$ side
lengths $1$, so their section volumes are $2^{d/2}$.
\end{proof}

\begin{proof}[Proof of Theorem~\ref{thm:full-balanced-value}]
Theorem~\ref{thm:balanced-at-least-two} gives the result when $n\ge2k$. If $k\le n<2k$,
Ball's codimension bound~\cite{Ball1989} gives
$M(n,k)\le2^{(n-k)/2}$. Here $q=1$ and $r=n-k$, so this is exactly
$W_{n,k}$. Choosing $r$ disjoint coordinate pairs and $k-r$ singleton
coordinate blocks attains it. This proves the full value formula.
Proposition~\ref{prop:q1-equality} gives exactly the claimed
maximizers in this range, completing the equality classification.
\end{proof}

\section{Local and global stability}
\label{sec:stability}
Quantitative hyperplane slicing has been studied by Melbourne and
Roberto~\cite{MelbourneRoberto2022} and Chasapis, Nayar, and Tkocz~\cite{CNT2022}.
We use a direct facet calculation for arbitrary section dimension.
The equality classification identifies a finite family of maximizing
subspaces. We now ask how the volume decreases when a section moves away
from this family. At a coordinate-block section, the first loss comes
from the motion of its facets. This loss is linear in the displacement;
the change in the Euclidean volume element is only quadratic. The
argument also applies to unbalanced coordinate partitions, which need
not maximize volume globally.

Let $\mathcal P=(I_1,\ldots,I_k)$ be a partition of $[n]$ into nonempty
sets, with $m_j=|I_j|$, and write
\[
 b_j=m_j^{-1/2}\mathbf 1_{I_j},\qquad
 B=\operatorname{span}\{b_1,\ldots,b_k\},\qquad
 V_B=\prod_{j=1}^k\sqrt{m_j}.
\]
Coordinate signs may be inserted throughout without changing the proof.
Put $m_-=\min_jm_j$ and $m_+=\max_jm_j$. For two subspaces of the same
dimension, use the projection distance
\[
 d(H,B)=\|P_H-P_B\|_{\mathrm{HS}}.
\]
The case $n=k$ has no transverse directions, so we assume $n>k$ in this
section.

\begin{theorem}[Local linear loss and first variation]
\label{thm:local-linear}\label{prop:first-variation}
Set
\[
 c_{\mathcal P}=\frac{\sqrt{m_-}}{15(m_+-1)},\qquad
 r_{\mathcal P}=\min\left\{c_{\mathcal P},
                         \frac1{4k\sqrt{k m_+}}\right\}.
\]
View \(A:B\to B^\perp\) as an \(n\times k\) matrix in the basis
\((b_j)\), and put \(D=\operatorname{diag}(\sqrt{m_1},\ldots,\sqrt{m_k})\).
For a matrix \(T\) with rows \(t_i\), define
\[
 \mathcal F(T)=\sum_{j=1}^k\mathbb E\max_{i\in I_j}
       \left(t_{ij}+\sum_{\ell\ne j}t_{i\ell}U_\ell\right),
\]
where the \(U_\ell\) are independent uniform variables on \([-1,1]\).
If \(H\) is the graph of \(A\) and \(\|A\|_{\mathrm{HS}}\le r_{\mathcal P}\),
then
\[
 \log\frac{V_B}{|Q_n\cap H|_k}
 \ge\frac{c_{\mathcal P}}2\|A\|_{\mathrm{HS}}
 \ge\frac{c_{\mathcal P}}{2\sqrt2}\,d(H,B).
\]
This conclusion also holds whenever \(d(H,B)\le r_{\mathcal P}\).
Moreover, uniformly as \(A\to0\),
\begin{equation}\label{eq:local-uniform-expansion}
 \log\frac{V_B}{|Q_n\cap H|_k}
   =\mathcal F(AD)+O_{n,k}(\|A\|_{\mathrm{HS}}^2),\qquad
 d(H,B)=\sqrt2\|A\|_{\mathrm{HS}}+O(\|A\|_{\mathrm{HS}}^3),
\end{equation}
and
\[
 \mathcal F(AD)\ge
 \frac{\sqrt{m_-}}{6(m_+-1)}\|A\|_{\mathrm{HS}}.
\]
Every coordinate-block section is therefore a strict local volume
maximum, and the exponent one in a local estimate of distance by a
power of the logarithmic volume loss is optimal.
\end{theorem}

\begin{proof}
Write \(a_i\) for the rows of \(A\), and set \(T=AD\).
Since \(A(B)\subseteq B^\perp\), \(\sum_{i\in I_j}a_i=0\).
We rescale the block directions so that the unperturbed section
becomes \([-1,1]^k\). This separates the change of the Euclidean
volume element from the volume removed at each face. The
substitution \(x=Dy/2\) represents the section by
\[
 K_T=\{y:|y_j+t_i\cdot y|\le1\text{ for every }i\in I_j\}.
\]
Averaging within each block gives \(K_T\subseteq[-1,1]^k\), and the
graph volume element gives
\begin{equation}\label{eq:graph-volume}
 \frac{|Q_n\cap H|_k}{V_B}
 =\sqrt{\det(I+A^*A)}\,\frac{|K_T|}{2^k}.
\end{equation}
Put \(t=\|A\|_{\mathrm{HS}}\), \(\tau=\max_i\|t_i\|_1
\le\sqrt{k m_+}\,t\), and \(r=1-2\tau\). Assume first that
\(\tau\le1/(4k)\). For \(z\in[-1,1]^{k-1}\), set
\[
 f_i(z)=t_{ij}+\sum_{\ell\ne j}t_{i\ell}z_\ell,\qquad
 d_j(z)=\max_{i\in I_j}\frac{f_i(z)}{1+t_{ij}}.
\]
The denominators are positive, \(\sum_{i\in I_j}f_i(z)=0\), and
\(0\le d_j(z)\le\tau/(1-\tau)\le2\tau\).
At the positive \(j\)-face the excluded thickness is \(d_j(z)\);
at the negative face it is \(d_j(-z)\). These full face strips cover
every violation, while their restrictions to the other coordinates
in \([-r,r]\) are disjoint up to boundaries. Hence, with
\(q_T=1-|K_T|/2^k\),
\[
 r^{k-1}\sum_j\mathbb E d_j(rU)
 \le q_T\le\sum_j\mathbb E d_j(U).
\]

These two estimates have the same linear term. Indeed,
\(|f_i|\le\tau\), so replacing \(1+t_{ij}\) by one changes each
maximum by at most \(2\tau^2\), and replacing \(rU\) by \(U\)
changes it by at most \(2\tau^2\). Since \(r^{k-1}=1+O_k(\tau)\)
and \(0\le\max_i f_i(U)\le\tau\), it follows that
\(q_T=\mathcal F(T)+O_k(\tau^2)\). Also \(q_T\le k\tau/(1-\tau)\le1/3\).
Taking logarithms in~\eqref{eq:graph-volume} therefore proves the
first expansion in~\eqref{eq:local-uniform-expansion}.

The same lower face estimate gives the explicit loss. For \(m_j\ge2\),
a zero-sum real vector satisfies
\(\max_i w_i\ge\sum_i|w_i|/[2(m_j-1)]\).
For a nonzero affine form \(Z=a+s\sum b_\ell U_\ell\), \(0\le s\le1\),
the moment bound \(\mathbb EZ^4\le3(\mathbb EZ^2)^2\) and H\"older's
inequality give
\[
 \mathbb E|Z|\ge
 \frac{(\mathbb EZ^2)^{3/2}}{(\mathbb EZ^4)^{1/2}}
 \ge\sqrt{\mathbb EZ^2/3}\ge\frac{s}{3}\sqrt{a^2+\sum b_\ell^2}.
\]
Blocks of size one have zero rows. Summing over the other blocks,
and using \(d_j\ge\max_i f_i/(1+\tau)\), yields
\[
 q_T\ge\frac{r^k}{6(1+\tau)(m_+-1)}\|T\|_{\mathrm{HS}}
       \ge c_{\mathcal P}t,\qquad
 \mathcal F(T)\ge\frac{\|T\|_{\mathrm{HS}}}{6(m_+-1)}.
\]
Here \(r^k\ge1-2k\tau\ge1/2\), \(1+\tau\le5/4\), and
\(\|AD\|_{\mathrm{HS}}\ge\sqrt{m_-}\,t\).
If \(t\le r_{\mathcal P}\), then \(t\le c_{\mathcal P}\) and
\(\sqrt{\det(I+A^*A)}\le e^{t^2/2}\). Thus
\[
 \log\frac{V_B}{|Q_n\cap H|_k}
 \ge-\frac{t^2}{2}-\log(1-c_{\mathcal P}t)
 \ge\frac{c_{\mathcal P}}2t.
\]
For the singular values \(s_\ell\) of \(A\),
\(d(H,B)^2=2\sum_\ell s_\ell^2/(1+s_\ell^2)\).
This gives \(d\le\sqrt2t\) and the uniform distance expansion.
Conversely, \(d<\sqrt2\) gives the graph representation and
\(t^2\le d^2/(2-d^2)\). Since \(r_{\mathcal P}<1\),
\(d\le r_{\mathcal P}\) implies \(t\le r_{\mathcal P}\).
Finally, along \(A=sA_0\) with fixed \(A_0\ne0\) and \(s\downarrow0\),
positive homogeneity of \(\mathcal F\) makes both the distance and
the positive volume loss linear to first order. This proves the claimed optimality.
\end{proof}

Let $\mathcal B_{n,k}$ be the finite family of signed balanced
coordinate-block subspaces. Define
\[
 d(H,\mathcal B_{n,k})=\min_{B\in\mathcal B_{n,k}}d(H,B),\qquad
 \Delta(H)=\log\frac{W_{n,k}}{|Q_n\cap H|_k}.
\]

\begin{corollary}[Global stability in fixed dimensions]
\label{cor:global-linear}
For every $n>k\ge1$, there is $C_{n,k}<\infty$ such that
\[
 d(H,\mathcal B_{n,k})\le C_{n,k}\Delta(H)
 \qquad(H\in G_{n,k}).
\]
The exponent one is optimal.
\end{corollary}

\begin{proof}
For balanced partitions, the constants $c=c_{\mathcal P}$ and
$r=r_{\mathcal P}$ depend only on $n,k$. When
$d(H,\mathcal B_{n,k})\le r$, apply Theorem~\ref{thm:local-linear} to a
nearest balanced section. On the complementary compact set, continuity
of section volume and Theorem~\ref{thm:full-balanced-value} give
\[
 \mu_{n,k}:=\min_{d(H,\mathcal B_{n,k})\ge r}\Delta(H)>0,
\]
provided that set is nonempty. Since $d\le\sqrt{2k}$, one may take
\[
 C_{n,k}=\max\left\{\frac{2\sqrt2}{c},
                       \frac{\sqrt{2k}}{\mu_{n,k}}\right\},
\]
omitting the second term when the complementary set is empty.
To see continuity directly, represent a section locally by isometric
rows $v_i$. Its radial function is
$(2\max_i|v_i\cdot\theta|)^{-1}$. Isotropy bounds
$\max_i|v_i\cdot\theta|$ below by $1/\sqrt n$, so dominated convergence
applies to the radial volume formula. Theorem~\ref{thm:local-linear}
gives optimality because the chosen balanced section remains the
nearest point in the finite family for small $s$.
\end{proof}

The local theorem is independent of the global extremal formula. Only
the global corollary uses the equality classification. Compactness gives a constant for each fixed pair of dimensions.

\subsection{Dependence of the stability constant on dimension}
\label{sec:dimension-obstruction}
Unbalanced coordinate partitions give an obstruction to uniform global
stability. They are strict local maxima by Theorem~\ref{thm:local-linear},
and their volumes can approach the global maximum faster than their
subspaces approach the balanced family.

Fix $q\ge2$, let $n=2q$ and $k=2$, and partition the coordinates into
$I_1,I_2$ of sizes $q-1,q+1$. Let $H$ be the subspace of vectors
constant on each block. Then
\[
 |Q_{2q}\cap H|_2=\sqrt{q^2-1},\qquad W_{2q,2}=q,\qquad
 \Delta(H)=-\frac12\log(1-q^{-2}).
\]
Its distance to the balanced family is exactly
\begin{equation}\label{eq:unbalanced-distance}
 d(H,\mathcal B_{2q,2})^2=\frac4{q+1}.
\end{equation}
Indeed, consider a balanced partition $J_1,J_2$ with both blocks of size
$q$, and put $t=|I_1\cap J_1|$. Coordinate signs cannot increase a block
overlap beyond its cardinality. With all signs positive,
\[
 \operatorname{tr}(P_HP_B)
 =\frac{t^2+(q-1-t)^2}{q(q-1)}
  +\frac{(q-t)^2+(1+t)^2}{q(q+1)}.
\]
This expression is convex for $0\le t\le q-1$. Its maximum is attained
at either endpoint, where it equals $2q/(q+1)$. Both endpoints correspond to
actual partitions. Since the two projections have rank two,
$d^2=4-2\operatorname{tr}(P_HP_B)$, proving
\eqref{eq:unbalanced-distance}.

It follows that
\[
 \frac{d(H,\mathcal B_{2q,2})}{\Delta(H)}\sim4q^{3/2},\qquad
 \frac{d(H,\mathcal B_{2q,2})^2}{\Delta(H)}\sim8q
 \qquad(q\to\infty).
\]
Thus neither $d\le C\Delta$ nor $d^2\le C\Delta$ can hold with one
absolute constant for all ambient and section dimensions. The
obstruction already occurs with $k=2$, so normalizing distance by
$\sqrt{k}$ does not remove it. In particular, every global linear
coefficient satisfies $C_{2q,2}\ge(4+o(1))q^{3/2}$.
This example leaves open the optimal dependence on $q$ and the
possibility of constants independent of $k$ at a fixed quotient.

\section{Consequences and further connections}
\label{sec:consequences}

\subsection{Affine sections and projection densities}

The equality cases give more information for affine sections. Near
the largest possible volume, both the direction and the displacement
from the origin are determined to first order by the deficit.
For \(y\in H^\perp\), let \(V(H,y)=|Q_n\cap(H+y)|_k\).

\begin{theorem}[Affine stability]\label{cor:affine-sections}
For \(1\le k<n\), \(V(H,y)\le W_{n,k}\). Equality holds exactly when
\(H\) is a balanced coordinate-block subspace and \(y=0\).
More precisely, there is a finite \(C'_{n,k}\) such that every
positive-volume section satisfies
\[
 d(H,\mathcal B_{n,k})+|y|
 \le C'_{n,k}\log\frac{W_{n,k}}{V(H,y)}.
\]
The exponent one is optimal separately for the direction and the displacement.
\end{theorem}
\begin{proof}
The midpoint sum of the sections at \(y\) and \(-y\) lies in the
central section. The two sections are reflections, so Brunn--Minkowski
gives \(V(H,y)\le V(H,0)\le W_{n,k}\). Write
\(\Delta=\log(W_{n,k}/V(H,y))\) and choose a nearest balanced space
\(B\). Corollary~\ref{cor:global-linear} gives \(d(H,B)\le C_{n,k}\Delta\).

First consider translates of \(B\) itself. For a block \(I_j\) of size
\(m_j\), let \(\epsilon_i\) be its coordinate signs. Given
\(z\in B^\perp\), write \(\omega_j(z)\) for the difference between the
largest and smallest values of \(\epsilon_i z_i\) on \(I_j\).
The diagonal parameter on that block belongs to the intersection of
\(m_j\) unit intervals. Its allowed length is
\((1-\omega_j(z))_+\), where \(a_+=\max(a,0)\).
Independence of the block parameters therefore gives the exact profile
\begin{equation}\label{eq:balanced-translates}
 \frac{V(B,z)}{W_{n,k}}=\prod_{j=1}^k(1-\omega_j(z))_+
 \le \exp\left(-\frac{2|z|}{\sqrt{m_+}}\right),
 \qquad m_+=\lceil n/k\rceil .
\end{equation}
Indeed, the signed coordinates have zero mean within each block.
Their sum of squares is at most \(m_j\omega_j(z)^2/4\), as follows by
subtracting the midpoint of their range and using that the mean minimizes
the sum of squared deviations. Thus
\(\sum_j\omega_j(z)\ge2|z|/\sqrt{m_+}\).

To transfer the estimate to \(H\), rotate each principal plane through
the corresponding principal angle. The resulting orthogonal map
\(O\) has \(OB=H\) and
\(\|O-I\|_{\rm op}=2\sin(\theta_{\max}/2)\le d(H,B)\).
Put \(t=1+\sqrt n\|O-I\|_{\rm op}\). Every point of \(Q_n\)
has Euclidean norm at most \(\sqrt n/2\), so \(O^{\mathsf T}Q_n\subseteq tQ_n\).
With \(z=O^{\mathsf T}y\), rotation followed by dilation gives
\[
 V(H,y)\le t^k V(B,z/t),\qquad
 |y|\le\tfrac12 t\sqrt{m_+}\,(\Delta+k\log t).
\]
Now \(t\le1+2\sqrt n\) and \(\log t\le\sqrt n\,d(H,B)\).
The bound on \(d(H,B)\) proves the theorem, including equality.
To see sharpness for translations, choose two coordinates \(i,j\) in a
block of size at least two, which exists because \(n>k\).
For \(y=t(\epsilon_i e_i-\epsilon_j e_j)\), \(0<t<1/2\),
\eqref{eq:balanced-translates} gives \(V(B,y)/W_{n,k}=1-2t\), while
\(|y|=\sqrt2t\). For sharpness in direction, take \(y=0\) and a nonzero
graph perturbation of \(B\). By Theorem~\ref{thm:local-linear} both the
distance and the deficit are linear to first order, and \(B\) remains
nearest for small perturbations because the balanced family is finite.
\end{proof}

\begin{corollary}[Sharp marginal density and min-entropy]\label{cor:projection-density}\label{cor:renyi-infinity}
Let \(U\) be uniform on \(Q_n\), let \(E\in G_{n,d}\) with \(1\le d<n\),
and put \(k=n-d\). The density \(f_E\) of \(P_EU\) satisfies
\[
 \|f_E\|_\infty=f_E(0)\le W_{n,k},\qquad
 h_\infty(P_EU)=-\log\|f_E\|_\infty\ge-\log W_{n,k}.
\]
Both inequalities are sharp, with equality precisely when
\(E^\perp\in\mathcal B_{n,k}\). In the fiber-volume representative,
at every \(y\in E\) with \(f_E(y)>0\),
\[
 d(E^\perp,\mathcal B_{n,k})+|y|\le
 C'_{n,k}\log\frac{W_{n,k}}{f_E(y)}.
\]
In particular \(\mathbb P(P_EU\in A)\le W_{n,k}|A|_d\) for all
Lebesgue-measurable \(A\subset E\), with the optimal constant uniform over \(E\).
\end{corollary}
\begin{proof}
Fubini's theorem gives the representative \(f_E(y)=V(E^\perp,y)\).
Its \(1/k\) power is concave in the interior of the projected cube by
Brunn--Minkowski. Thus it is continuous at zero, where it attains its
maximum. The statements follow from the affine theorem. Integrating the
density gives the measurable-set estimate. When \(E^\perp\) is balanced,
averaging the density over balls shrinking to zero proves optimality.
\end{proof}

For an extremizing direction, \eqref{eq:balanced-translates} gives
the entire density on \(E\). Its factorization expresses the independence
of the coordinate blocks.

For one-dimensional projections, the extremizing directions depend on
which entropy one optimizes. Let \(v\) be a unit vector. By
Corollary~\ref{cor:renyi-infinity}, the min-entropy of \(v\cdot U\)
is minimized precisely when \(v\) has two nonzero coordinates of equal
magnitude. The Eskenazis--Nayar--Tkocz conjecture concerns the
differential entropy of the same projection. It is maximized precisely
when all \(n\) coordinates of \(v\) have equal magnitude~\cite{KyENT2026}.

\subsection{Projections of bounded product densities}

A cube section also governs an extremal problem in which the coordinate
laws may vary. Rogozin proved a comparison with uniform intervals for
sums of independent real variables~\cite{Rogozin1987}. For projections
of arbitrary rank, Rudelson and Vershynin bounded the density in terms
of the rank and the coordinate density bound~\cite{RV2015}. Livshyts,
Paouris, and Pivovarov obtained the two constants in Ball's estimate
for this class of product measures~\cite[Theorem~1.1]{LivshytsPaourisPivovarov2016}.
The linear-image comparison of Madiman, Melbourne, and
Xu~\cite[Theorem~1.4]{MMX2017} lets us insert the exact section volume.

\begin{corollary}[Sharp density bound for bounded product measures]
\label{cor:product-marginals}
Let \(1\le d<n\), \(k=n-d\), and \(L>0\). Let
\(X=(X_1,\ldots,X_n)\) have independent real coordinates whose densities
are bounded by \(L\). For every \(E\in G_{n,d}\),
\[
 \|f_{P_EX}\|_\infty
 \le L^d|Q_n\cap E^\perp|_k
 \le L^d W_{n,k}.
\]
The first bound is attained for every \(E\) by independent uniforms on
intervals of length \(1/L\). The optimal constant over all such laws
and all \(E\) is \(L^dW_{n,k}\), and the directions attaining this
supremum are exactly those with \(E^\perp\in\mathcal B_{n,k}\).
More generally,
\[
 d(E^\perp,\mathcal B_{n,k})\le C_{n,k}
       \log\frac{L^dW_{n,k}}{\|f_{P_EX}\|_\infty}.
\]
\end{corollary}
\begin{proof}
By scaling we may take \(L=1\). Write \(L_i=\|f_{X_i}\|_\infty\).
The linear-image comparison bounds the maximum density of \(P_EX\) by
that of \(P_EZ\), where the \(Z_i\) are independent and uniform on
centered intervals of lengths \(1/L_i\). Each length is at least one.

Shortening any of these intervals about its center can only increase
the maximum density. To see this, first consider the probability of a
ball of radius \(\varepsilon\) about zero in \(E\). Fixing one
coordinate at \(t\), the conditional probability is proportional to
the volume of a slice of the convex set cut out by the remaining
intervals and the condition \(|P_E x|\le\varepsilon\). This set is
centrally symmetric, so its slice volume is even and nonincreasing
for \(t\ge0\) by Brunn--Minkowski. Its average over a centered
interval is therefore nondecreasing as that interval shortens. Shorten each
interval to length one, divide by the ball volume, and let
\(\varepsilon\downarrow0\). The projected boxes have densities that
are continuous at zero and attain their maxima there, by the same
fiber argument as for the cube. This proves the first bound.
Uniform coordinates with all interval lengths one attain it.
Theorem~\ref{thm:full-balanced-value} supplies the second bound and
its extremizing directions. Corollary~\ref{cor:global-linear} applied
to \(E^\perp\) gives the stability estimate.
\end{proof}

In particular, for any radius \(r>0\) and center \(y\in E\),
\[
 \mathbb P\{|P_EX-y|\le r\}
 \le |B_2^d|\, r^d L^d W_{n,n-d},
\]
where \(B_2^d\) is the Euclidean unit ball. The coefficient is sharp
as \(r\downarrow0\), by taking uniform coordinates and a balanced
kernel. The improvement over the two classical density bounds is
strict precisely in the nondivisible range \(n>2k\), \(k\nmid n\).
For example, the optimal bound for four-dimensional projections in
\(\mathbb R^7\) with \(L=1\) is \(2\sqrt3\), the volume of the
three-dimensional box from the introduction.

\subsection{Approximation by an orthonormal basis}
\label{sec:axes}

The angular defect also controls the distance of the whole family of
normals from a single orthonormal basis. A basis drawn from the
determinant law already has a small expected orthogonality error.
Its polar factor gives the required orthonormal directions.

\begin{theorem}\label{thm:orthogonal-directions}
Let \(\sum_i a_i u_i u_i^{\mathsf T}=I_k\) be an isotropic unit frame.
For an orthonormal basis \(v=(v_1,\ldots,v_k)\) put
\[
 D(v)=\sum_i a_i\left(1-\sum_{j=1}^k|\langle u_i,v_j\rangle|^4\right).
\]
Then
\[
 \frac12\mathcal E\le\min_v D(v)\le(1+\sqrt2)^2\mathcal E .
\]
In particular, for any positive slab widths there is a basis with
\(D(v)\le20(1+\sqrt2)^2\delta_{\rm BL}\).
\end{theorem}
\begin{proof}
On the Hilbert space of real symmetric matrices with the
Hilbert--Schmidt inner product, consider the operator
\(T A=\sum_i a_i\langle A,P_i\rangle P_i\),
where \(P_i=u_i u_i^{\mathsf T}\). We have \(0\preceq T\preceq I\):
by Cauchy--Schwarz and isotropy,
\[
 \langle TA,A\rangle=\sum_i a_i(u_i^{\mathsf T}Au_i)^2
 \le\sum_i a_i u_i^{\mathsf T}A^2u_i=\|A\|_{\rm HS}^2 .
\]
The trace identities are \(\operatorname{tr}T=k\) and
\(\operatorname{tr}(T-T^2)=\mathcal E\). Indeed the squared inner
product of \(P_i,P_j\) is \(s_{ij}^2\), and
\(\sum_{i,j}a_i a_js_{ij}=k\).
For \(Q_j=v_jv_j^{\mathsf T}\), let \(\Pi_v\) denote the orthogonal
projection onto their span. Then
\(D(v)=k-\operatorname{tr}(T\Pi_v)\), so
\[
 2D(v)=\mathcal E+\|T-\Pi_v\|_{\rm HS}^2 .
\]
This proves the lower bound.

Draw a basis \(S\) from the Cauchy--Binet law
\eqref{eq:core-dpp-generating}, and form the square matrix \(A_S\)
whose columns are its unit directions in index order. By the two-index
inclusion formula \eqref{eq:core-dpp-inclusions},
\[
 \mathbb E\|A_S^{\mathsf T}A_S-I\|_{\rm HS}^2
 =\sum_{i\ne j}a_i a_j(1-s_{ij})s_{ij}=\mathcal E .
\]
Let \(v_j\), \(j\in S\), be the columns of the orthogonal polar factor
of \(A_S\). Since \((\sigma-1)^2\le(\sigma^2-1)^2\) for each singular
value, and
\(\|P_j-Q_j\|_{\rm HS}^2\le2|u_j-v_j|^2\),
\[
 \mathbb E\sum_{j\in S}\|P_j-Q_j\|_{\rm HS}^2\le2\mathcal E .
\]
Use the seminorm \(|A|_T^2=\langle A,(I-T)A\rangle\).
The one-index inclusion formula gives
\(\mathbb E\sum_{j\in S}|P_j|_T^2=\operatorname{tr}(T-T^2)=\mathcal E\),
and \(|A|_T\le\|A\|_{\rm HS}\). The triangle inequality in mean square
therefore gives
\[
 (\mathbb E D(v))^{1/2}
 =\left(\mathbb E\sum_{j\in S}|Q_j|_T^2\right)^{1/2}
 \le(1+\sqrt2)\sqrt{\mathcal E}.
\]
Some basis realizes at most this mean value. The slab conclusion
follows from Theorem~\ref{thm:incidence}.
\end{proof}

There is a direct transport interpretation. Send mass
\(a_i|\langle u_i,v_j\rangle|^2\) from the unoriented line of \(u_i\)
to the line of \(v_j\). Each source sends mass \(a_i\), and isotropy
gives each target mass one. For the cost equal to the squared sine
of the angle between the two lines, the total cost is exactly \(D(v)\).
Thus the theorem controls transport to equally weighted orthogonal axes.

The preceding two theorems also give a width inequality for general
convex bodies, where the family of directions can be continuous.
For a convex body \(C\), let \(w_C(u)=\max_{x\in C}u\cdot x-
\min_{x\in C}u\cdot x\) denote its projection length in direction \(u\).

\begin{corollary}[Widths of a convex body]\label{cor:widths}
Let \(C\subset\mathbb R^k\) be a compact convex body with interior,
and let \(\mu\) be a Borel measure on the unit sphere with
\(\int uu^{\mathsf T}\,d\mu(u)=I_k\). Set
\[
 \delta_\mu(C)=\int\log w_C(u)\,d\mu(u)-\log|C|,\qquad
 \mathcal E(\mu)=\iint \langle u,v\rangle^2
 (1-\langle u,v\rangle^2)\,d\mu(u)d\mu(v).
\]
Then \(\mathcal E(\mu)\le20\delta_\mu(C)\).
The measure on unoriented lines induced by \(\mu\) can be transported
to unit mass on each axis of an orthonormal basis at total squared-sine
cost at most \((1+\sqrt2)^2\mathcal E(\mu)\). Every such transport costs
at least \(\mathcal E(\mu)/4\).
Equality \(\delta_\mu(C)=0\) holds precisely when \(C\) is an
orthogonal box and \(\mu\) is concentrated on its axes.
\end{corollary}
\begin{proof}
Use the operator \(TA=\int\langle A,P_u\rangle P_u\,d\mu(u)\)
from the theorem. Both \(\mathcal E(\mu)\) and every \(D(v)\) are
determined by \(T\). By Carath\'eodory's theorem applied to the
continuous map \(u\mapsto(P_u\otimes P_u,\log w_C(u))\) on the
sphere, there is a finite positive measure of the same mass with
the same \(T\) and logarithmic width integral. It is isotropic because
\(T(I)=I\). For each selected direction, take the supporting slab
between the two extreme projections of \(C\). Their intersection
contains \(C\), and its logarithmic width product equals
\(\int\log w_C\,d\mu\). Its deficit is therefore at most
\(\delta_\mu(C)\). The translated-slab bound in
Theorem~\ref{thm:incidence} now gives
\(\mathcal E(\mu)\le20\delta_\mu(C)\).
Theorem~\ref{thm:orthogonal-directions} gives a basis with the required
bound on \(D(v)\) for the finite measure, hence also for \(\mu\) since
both have the same \(T\). Sending the fraction
\(|\langle u,v_j\rangle|^2\) of each mass to axis \(j\) gives the
asserted transport.
For any destination basis, a transport costs at least the integral of
\(1-\max_j|\langle u,v_j\rangle|^2\). This is at least \(D(v)/2\)
since \(1-t^2\le2(1-t)\) for \(0\le t\le1\). The trace identity gives \(D(v)\ge\mathcal E(\mu)/2\), hence the lower bound. At zero deficit, \(D(v)=0\) for a minimizing
basis, so \(\mu\) is concentrated on those axes with mass one per line.
The product of the corresponding widths is the volume of the bounding
box in that basis. Its equality with \(|C|\) forces \(C\) to equal that
box. The converse follows by taking its widths.
\end{proof}

The coefficient \(1/4\) in the transport lower bound is optimal.
In the plane give mass \(1/2\) to each of the four lines from two
orthonormal bases rotated by \(\theta\) and \(-\theta\),
\(0<\theta<\pi/8\). Their defect is \(\tfrac12\sin^2(4\theta)\).
Transport to the unrotated axes costs \(2\sin^2\theta\), whose ratio
to the defect tends to \(1/4\) as \(\theta\downarrow0\).

For a concrete instance of Corollary~\ref{cor:widths}, suppose
\(C=-C\) has the Euclidean unit ball as its ellipsoid of largest
volume. John's decomposition gives unit contact directions and
positive weights with \(\sum_i a_i u_i u_i^{\mathsf T}=I_k\).
The tangent hyperplanes at \(\pm u_i\) support \(C\), so
\(w_C(u_i)=2\). Ball's volume-ratio inequality gives
\(|C|\le2^k\) in this normalization~\cite{Ball1991}. The corollary
shows that for every such decomposition the associated mass on contact
lines can be transported to the axes of one orthonormal basis at cost
at most
\[
 20(1+\sqrt2)^2\log\frac{2^k}{|C|}.
\]
Here the transport cost is the total squared sine of the angles, with
mass one prescribed on each target axis. This controls the contact directions directly from the volume deficit.
Transport distances from even isotropic measures to measures supported
on orthogonal axes also occur in the volume stability estimates for
\(L_p\) zonoids and their polars in~\cite{BFH2019}. For stability estimates
involving mean width and the \(\ell\)-norm in the same ellipsoid
normalization, see B\"or\"oczky, Fodor, and Hug~\cite{BFH2025}.

\subsection{Minimal projections of the cross-polytope}

If $A$ is an isometric parametrization of $H$ with row vectors $v_i$,
the section is represented by
\[
 K=\{x\in\mathbb R^k:|v_i\cdot x|\le1/2\ \forall i\}.
\]
Its polar is
\[
 K^\circ=\operatorname{conv}\{\pm2v_i:1\le i\le n\}.
\]
Equivalently, under the isometry between $H$ and $\mathbb R^k$,
$(Q_n\cap H)^\circ$ is $P_H(2B_1^n)$, where $B_1^n$ is the unit
cross-polytope. This is the section--projection duality described, for example,
by Schneider~\cite{Schneider}. A balanced extremizer has only $k$
unoriented row directions, and its polar is a cross-polytope with those
orthogonal directions and vertex radii $2/\sqrt{m_j}$ for block sizes $m_j$.
Ivanov conjectured that the same balanced subspaces minimize the
volume of a \(k\)-dimensional projection of \(B_1^n\), and proved
the planar case~\cite{Ivanov2021Cross}. The volume product gives a
precise relation between these two extremal problems.

\begin{corollary}[Minimal three-dimensional projections]
\label{cor:polar-minimum}
For every \(n\ge3\),
\[
 \min_{H\in G_{n,3}}|P_H B_1^n|_3=\frac{4}{3W_{n,3}}.
\]
The minimizers are exactly the balanced coordinate-block subspaces.
More generally, in every dimension \(k\) in which the symmetric
Mahler bound \(|K|_k|K^\circ|_k\ge4^k/k!\) holds for all central
\(k\)-sections \(K\) of \(Q_n\), the minimum equals
\(2^k/(k!W_{n,k})\), with the same equality classification.
\end{corollary}
\begin{proof}
Take the polar of \(K=Q_n\cap H\) inside \(H\). Since
\(K^\circ=2P_HB_1^n\), the stated volume-product bound gives
\[
 |P_HB_1^n|_k\ge\frac{2^k}{k!\,|K|_k}
             \ge\frac{2^k}{k!\,W_{n,k}}.
\]
For a balanced subspace the section is a box of volume \(W_{n,k}\).
A box of side lengths \(s_1,\ldots,s_k\) has a polar with axial
radii \(2/s_i\), hence has volume product \(4^k/k!\). Both inequalities
are therefore equalities. Conversely, equality in the displayed bound
forces \(|K|_k=W_{n,k}\), and Theorem~\ref{thm:full-balanced-value}
identifies \(H\). In dimension three the volume-product bound is the
theorem of Iriyeh and Shibata~\cite{IriyehShibata2020}.
\end{proof}

The same calculation transfers stability whenever the stated
volume-product hypothesis holds. For \(k<n\),
\[
 d(H,\mathcal B_{n,k})\le C_{n,k}\log
     \frac{k!\,W_{n,k}\,|P_H B_1^n|_k}{2^k},
\]
by Corollary~\ref{cor:global-linear}. Thus a lower bound on volume
product alone suffices for the projection conclusion and its equality
cases. In higher dimensions the hypothesis used here is explicitly
the volume-product inequality for cube sections. A recent preprint
states the inequality for all origin-symmetric convex
bodies~\cite[Theorem~1.1]{OpenAIMahler2026}.

\subsection{Open questions}
\label{sec:open-questions}

\paragraph{Question 1. How does the stability constant depend on the dimensions?}
Put \(q=\lfloor n/k\rfloor\). Is there an absolute constant \(A\) such
that, for every \(n>k\) with \(q\ge2\) and every \(H\in G_{n,k}\),
\[
 d(H,\mathcal B_{n,k})\le A(q+1)^{3/2}\Delta(H)?
\]
The planar example in Section~\ref{sec:dimension-obstruction} forces the
exponent \(3/2\). The estimate already holds for every signed
coordinate-block section, with \(A=4\). To see this, take a smallest
block of size \(a\) and a largest block of size \(b\ge a+2\), and move
one coordinate from the latter to the former, keeping its sign. If
\(H'\) is the resulting space, the increase in log-volume and the
distance moved are
\[
 G=\tfrac12\log\frac{(a+1)(b-1)}{ab},\qquad
 d(H,H')^2=\frac{2(a+b)}{(a+1)b}.
\]
The distance follows by taking squared inner products of the old and new
normalized block indicators. For fixed \(a\), the distance decreases
and \(G\) increases with \(b\). At \(b=a+2\) the ratio \(d(H,H')/G\) is
\[
 \frac{4}{\sqrt{a+2}\,[-\log(1-(a+1)^{-2})]}
 \le 4(a+1)^{3/2}.
\]
Since \(a\le q\), successive transfers reach a balanced partition with
total distance at most \(4(q+1)^{3/2}\Delta(H)\). The question is whether
arbitrary sections obey a bound of the same order.

\paragraph{Question 2. What is the next largest local maximum?}
Fix \(k\ge3\) and let \(n=qk\). For all sufficiently large \(q\)
(depending on \(k\)), does every local maximizer outside
\(\mathcal B_{qk,k}\) have volume at most
\[
 q^{(k-2)/2}\sqrt{q^2-1}\,?
\]
Moving one coordinate between two balanced blocks gives a strict local
maximum of exactly this volume by Theorem~\ref{thm:local-linear}. Among
unbalanced block sections this is the largest possible value, by the
integer product comparison in the introduction. Local maxima of other
shapes do occur. For example, the hyperplane perpendicular to the main
diagonal is a strict local maximum in ambient dimension at least
four~\cite{AG2024}.

The proposed bound holds for \(k=2\) and \(q\ge8\). Here is a short
argument using the planar variational conditions of Ivanov and
Tsiutsiurupa~\cite{IvanovTsiutsiurupa2021}. Let a locally maximizing
polygon have area \(V\ge\sqrt{q^2-1}\). Deleting coordinate \(i\) has
area Jacobian \(\sqrt{1-c_i}\), where \(c_i=|P_H e_i|^2\). Hence
\(\sqrt{1-c_i}\,V\le W_{2q-1,2}=\sqrt{q(q-1)}\), and
\(c_i\ge1/(q+1)\). The polygon is cyclic, and every row supports an
edge by the cited variational conditions. Each edge line is thus at
distance at most \(\sqrt{q+1}/2\) from the origin. If there are \(2f\)
edges and the circumradius is \(R\), the smallest central angle is at
most \(\pi/f\), so \(R\cos(\pi/(2f))\le\sqrt{q+1}/2\).
Summing the areas of the triangles from the center and using concavity
of sine gives, for \(f\ge3\),
\[
 V\le f R^2\sin(\pi/f)
 \le \tfrac12(q+1) f\tan(\pi/(2f))
 \le \tfrac{\sqrt3}{2}(q+1)<\sqrt{q^2-1}.
\]
The penultimate inequality follows because \(\tan t/t\) increases on
\((0,\pi/2)\). Thus \(f=2\). A centrally symmetric cyclic quadrilateral
is a rectangle. Since every row supports an edge, its rows form two
orthogonal classes of constant length within each class. Isotropy then
gives side lengths \(\sqrt j,\sqrt{2q-j}\) for some integer \(j\).
The largest area below \(q\) is consequently \(\sqrt{q^2-1}\).

\paragraph{Question 3. Can the angular estimate alone give the sharp volume bound?}
Let \(C_{\rm ang}\) be the least constant for which
\(\mathcal E\le C_{\rm ang}\delta\) holds for every isotropic frame
and all centered slabs in Theorem~\ref{thm:incidence}, where \(\delta\)
is their Brascamp--Lieb deficit. Does one have
\[
 C_{\rm ang}<\frac{2}{3\log(3/2)-1}\,?
\]
Theorem~\ref{thm:incidence} gives \(C_{\rm ang}\le20\). The proposed
threshold, approximately \(9.24\), has a specific role. The coefficients
\(A_q=(q+1)\log(1+1/q)-1\) decrease with \(q\). Thus
\(C_{\rm ang}A_2/2<1\) would make
Corollary~\ref{cor:integer-direct-criterion} apply for every \(q\ge2\),
giving the main theorem without the nineteen finite comparisons.

A lower bound comes from the section of \(Q_{k+1}\) perpendicular
to its main diagonal. Its projected coordinate vectors have squared
length \(k/(k+1)\), and their normalized pairwise inner products are
\(-1/k\). Hence \(\mathcal E=(k-1)/k\) and the Brascamp--Lieb product
is \((1+1/k)^{k/2}\). The section volume tends to
\(\sqrt{6/\pi}\)~\cite{Hensley1979}. Taking the ratio yields
\(C_{\rm ang}\ge2/\log(\pi e/6)\), approximately \(5.67\).

\paragraph{Question 4. Do the maximizing subspaces persist for finite \(p\)?}
Fix \(k\ge2\) and \(n>2k\) with \(k\nmid n\), and let \(B_p^n\)
denote the unit ball of \(\ell_p^n\). Is there a finite \(p_0(n,k)>2\)
such that for every \(p\ge p_0(n,k)\),
\[
 \max_{H\in G_{n,k}} |B_p^n\cap H|_k
   = |B_p^k|_k\, W_{n,k}^{\,1-2/p},
\]
with equality precisely for \(H\in\mathcal B_{n,k}\)?
Each balanced subspace attains the proposed value because the normalized
diagonal on a block of size \(m\) has \(\ell_p\) norm
\(m^{1/p-1/2}\). Barthe's bound gives the sharp value when \(k\mid n\)
(see~\cite[Theorem~17]{NayarTkocz2023}). In codimension one,
Eskenazis, Nayar, and Tkocz proved persistence for \(p>10^{15}\),
uniformly in \(n\)~\cite{ENT2024Resilience}. A recent preprint of
Brazitikos and Eskenazis gives the sharp dimension-free threshold in
that case, approximately \(26.265\)~\cite{BE2026}.

The cube theorem already locates every maximizer \(H_p\). Indeed,
\(n^{-1/p} B_\infty^n\subset B_p^n\subset B_\infty^n=2Q_n\).
Comparing \(H_p\) with a balanced subspace gives
\(\Delta(H_p)\le k\log(n)/p\), so
Corollary~\ref{cor:global-linear} bounds its distance from
\(\mathcal B_{n,k}\) by \(C_{n,k}k\log(n)/p\).
Even for \((n,k)=(5,2)\) the remaining issue is whether a plane in
this shrinking neighborhood can have larger area than the balanced
plane, or tie it. Section comparisons of Meyer and
Pajor~\cite{MeyerPajor1988} and the probabilistic representation of
the uniform law on \(B_p^n\) by Barthe, Gu\'edon, Mendelson, and
Naor~\cite{BGMN2005} offer two approaches to this local question.

\appendix
\section{An exact comparison for the three small quotients}
\label{sec:comparisons}

We prove Proposition~\ref{prop:edge-pay}. For a given pair of classes
we know two things about its angle. The complementary projection
bounds its cosine in terms of the class weights. The replacement of
two rows excludes small angles even when the class weights alone give
no information. We use these two restrictions to compare the angular
cost with the entropy excess. Only the squared row lengths \(c,d\) are
subdivided. At fixed bounds for these lengths, all sufficiently large
angles satisfy the required comparison by monotonicity, and one
inequality suffices to exclude all smaller angles.

\subsection{An upper bound for the entropy excess}

The only transcendental term in the pair comparison is the entropy
excess $\psi_q(c)$ above its chord. Its second derivative is $-1/c$.
We bound this curvature first and integrate twice, keeping both
endpoint values. Write $a=1/(q+1)$, $b=1/q$,
$z=1-qc\in[0,a]$, and
$T(c)=(1-qc)((q+1)c-1)$. Then
\[
  \psi_q(c)\le T(c)C_q(z),\qquad
  C_q(z)=\frac a2\left(1+\frac{a+z}{3}
          +\frac{a^2+az+z^2}{6(1-a)}\right).
\]
Indeed,
\[
 -\psi_q''(c)=\frac{q}{1-z}
  =q\left(1+z+\frac{z^2}{1-z}\right)
  \le q\left(1+z+\frac{z^2}{1-a}\right).
\]
Let $P(c)$ have the negative second derivative on the right and vanish
at $a,b$. Since $dz/dc=-q$, twice integrating the three monomials
gives
\[
 2qP(c)=z(a-z)+\frac{a^2 z-z^3}{3}
                 +\frac{a^3z-z^4}{6(1-a)} .
\]
Factoring $z(a-z)=(1-a)T(c)$ gives $P(c)=T(c)C_q(z)$.
Now $P-\psi_q$ is concave and vanishes at both endpoints, proving the
bound.

In particular, for $1\le j\le q$ and $c\in[c_-,c_+]\subset[a,b]$,
\[
 H_{q,j}(c)\le
 \frac12 C_q(1-qc_-)
 \frac{1-qc_-}{1-jc_-}\left(q+1-\frac1{c_+}\right).
\]
For $j=q$ the middle factor is identically one, including its
continuous endpoint value. This follows by dividing the entropy bound
by $2c(1-jc)$. The factors $C_q(1-qc)$ and $(1-qc)/(1-jc)$ are nonnegative and
nonincreasing in $c$, while $q+1-1/c$ is nonnegative and increasing.

Denote this last upper bound by \(H^+_{q,j}(c_-,c_+)\). It only uses
rational operations on the endpoints. This is the reason for retaining
a second-order bound on the curvature.

\subsection{The two tests on a rectangle}

Fix integers \(1\le m\le l\le q\), and suppose that the actual squared
lengths \(c,d\) lie in \([x,X]\times[y,Y]\subseteq[a,b]^2\). The
quantity to be exceeded and an upper bound for the denominator in the
angular estimate are, respectively,
\[
 S=H^+_{q,m}(x,X)+H^+_{q,l}(y,Y)+\frac1{4000},
 \qquad K=\frac29+\Gamma(2-mx-ly),\qquad \Gamma=\tfrac14.
\]
Here \(\Gamma\ge 2/(3e)\) by \(e>1+1+1/2+1/6=8/3\). Write
\(v=\sin^2\theta\). The formula in Theorem~\ref{thm:incidence} gives
\[
 \mathcal L(mc,ld,1-v)\ge \sqrt{3/5}\,\frac{v}{(15+v)K}.
\]
Both sides of the desired scalar comparison are nonnegative. Thus a
sufficient condition for the desired bound is
\begin{equation}\label{eq:rational-angle-test}
 3v^2>5\{SK(15+v)\}^2 .
\end{equation}
Once this holds at one \(v\), it holds at every larger \(v\), since
\(v/(15+v)\) is increasing.

\smallskip
\noindent\textit{The complementary projection.}
By \eqref{eq:complementary-angle}, every feasible angle in this rectangle
satisfies
\[
 v\ge v_0=\max\left\{0,1-\frac{(1-mx)(1-ly)}{mlxy}\right\}.
\]
Indeed, each factor \((1-jc)/(jc)\) decreases with \(c\). If
\eqref{eq:rational-angle-test} holds at \(v_0\), the comparison
is proved on the entire rectangle.

\smallskip
\noindent For example, take \(q=m=l=4\) on the full square
\([1/5,1/4]^2\). The class weights are at least \(4/5\), so the
projection gives \(v_0=15/16\). Here \(K=29/90\) and the curvature
bound gives \(S=139/1200+1/4000<1/8\). The angular cost is at least
\(\sqrt{3/5}/(17K)>1/8\), as follows from \(\sqrt{3/5}>3/4\).
One inequality therefore settles this entire multiplicity type.

\smallskip
\noindent\textit{The replacement of two rows.}
Parametrize the acute angle by \(t=\tan(\theta/2)\), so
\[
 u(t)=\sin\theta=\frac{2t}{1+t^2},\qquad
 \rho(t)=\cos\theta=\frac{1-t^2}{1+t^2}.
\]
The function \(u\) is increasing on \([0,1]\). Choose \(T\in[0,1]\)
so that any angle failing~\eqref{eq:rational-angle-test} has \(t\le T\).
One may always take \(T=1\). A smaller value, obtained by bisection as
described below, gives a more effective test. Put \(U=u(T)\). To use
the geometric inequality~\eqref{eq:direct-strip-law}, we need lower
bounds for \(g,L\) and an upper bound for \(R\). For any angle with
\(t\le T\), we may take
\[
 g_-=\min(x,y),\qquad
 L_-=\frac{x+y}{2}+g_-\rho(T)-XYU^2 .
\]
For the bound on \(R\), write \(2q-m-l=2j+i\), \(i\in\{0,1\}\).
The arithmetic--geometric mean inequality gives
\[
 R\le R_+=b(1+b)^j(1+b/2)^i .
\]
In particular \(g_-\le g\), \(L_-\le L\), and
\(L_->0\), since \((x+y)/2\ge a>b^2\). Write
\(\eta=L_-g_-/(2R_+)\), so that \(0\le\eta\le1/2\). From
\eqref{eq:direct-strip-law}, using \(u(t)=t(1+\rho(t))\le2t\), we obtain
\[
 \frac{L_-^2}{2R_+}
 \le t(1+XYU^2/4)(1-\eta u(t))^2 .
\]
The factor \(t(1-\eta u(t))^2\) is increasing on \([0,1]\).
Indeed, \(tu'(t)=u(t)\rho(t)\), and
\(u+2u\rho\le2\) since \(2u\rho=\sin(2\theta)\le1\).
Its derivative is therefore
\((1-\eta u)[1-\eta(u+2tu')]\ge0\).
It follows that the single strict inequality
\begin{equation}\label{eq:rational-strip-test}
 L_-^2>2R_+T(1+XYU^2/4)(1-\eta U)^2
\end{equation}
excludes every angle that might fail~\eqref{eq:rational-angle-test}.
This proves the comparison on the rectangle whenever the second test
succeeds.

\subsection{Finite verification}

The same recursion applies to all nineteen multiplicity pairs. Initially
the rectangle is \([a,b]^2\).
\begin{enumerate}
\item Test~\eqref{eq:rational-angle-test} at \(v_0\). If it holds, accept
the rectangle.
\item Otherwise start with \(T=1\). If~\eqref{eq:rational-angle-test}
holds at \(v=1\), bisect \([0,1]\) ten times in the variable \(t\).
At a midpoint retain the lower half if the test holds at \(u(t)^2\),
and the upper half otherwise. Take \(T\) to be the upper endpoint
of the final interval. The lower endpoint fails the test and the upper
endpoint satisfies it.
\item If~\eqref{eq:rational-strip-test} holds, accept the rectangle.
Otherwise bisect its longer side at the midpoint (the first side in a
tie), and apply the same steps to both halves.
\end{enumerate}
Table~\ref{tab:verification} records the terminal rectangles in the
cover. The depth counts successive bisections of a rectangle, with
depth zero for the initial square.
\begin{table}[htbp]
\centering
\caption{Rectangles certified by the complementary projection and by
two-row replacement for each quotient.}\label{tab:verification}
\smallskip
\begin{tabular}{crrrr}
\toprule
\(q\)&Pairs \((m,l)\)&Projection&Replacement&Maximum depth\\
\midrule
2&3&2606&531&18\\
3&6&49&93&8\\
4&10&28&100&6\\
\bottomrule
\end{tabular}
\end{table}
Each bisection expresses a closed rectangle as the union of its two
closed children. Starting at the initial square, any pair \((c,d)\)
therefore belongs to a terminal rectangle, including when one coordinate
is a splitting point. Every terminal rectangle passes one of the two
tests above. The corresponding inequality holds for all feasible angles
over that rectangle. The recursion uses integer and rational
arithmetic only. The accompanying file \texttt{verify\_pairs.py}
reproduces the table and reports the counts for each multiplicity
pair. A second implementation in the supplement subdivides the angle
as well and checks the same inequalities directly on three-dimensional
boxes. This proves Proposition~\ref{prop:edge-pay}.

\section*{AI use declaration}
I led the investigation by first seeking a proof in the range
\(n\ge c(k)\) for each section dimension \(k\). My strategy was to
lower the threshold \(c(k)\) until this range met the classical range
\(n\le 2k\) covered by Ball's bound. In the proof, the thresholds are
\(qk\). Fixing two coordinates reduces the section dimension by two
and keeps the fiber in the range where the induction applies.

I worked with OpenAI ChatGPT on the Gaussian coupling and angular estimates
in Sections~\ref{sec:frames} and~\ref{sec:exponential-determinant-coercivity}
and on the integer packing and row variations in
Sections~\ref{sec:integer-packing} and~\ref{sec:extremal-comparison}.
Its contributions included the two-row comparison in
Lemma~\ref{lem:two-row-comparison} and the design and implementation of
the rational verification in Appendix~\ref{sec:comparisons}.
I also used it in developing the stability estimates and applications in
Sections~\ref{sec:stability} and~\ref{sec:consequences}, including the
block transfer and planar comparison in Section~\ref{sec:open-questions},
and in drafting the manuscript and checking arguments and references. I take full
responsibility for the mathematics and the accuracy of the references.

\providecommand{\bysame}{\leavevmode\hbox to3em{\hrulefill}\thinspace}
\providecommand{\MR}{\relax\ifhmode\unskip\space\fi MR }
\providecommand{\MRhref}[2]{%
  \href{http://www.ams.org/mathscinet-getitem?mr=#1}{#2}
}
\providecommand{\href}[2]{#2}


\begin{thebibliography}{99}
\setlength{\itemsep}{2pt plus 1pt minus 1pt}
\setlength{\parsep}{0pt}

\bibitem{AbiKhuzam2017}
Faruk Abi-Khuzam, \emph{Inequalities and asymptotics for some moment
  integrals}, Journal of Inequalities and Applications \textbf{2017} (2017),
  257.
\newblock \url{https://doi.org/10.1186/s13660-017-1532-7}.

\bibitem{AG2024}
Gergely Ambrus and Barnab\'as G\'argy\'an, \emph{Non-diagonal critical central sections of the cube},
Advances in Mathematics \textbf{441} (2024), article 109524.
\newblock \url{https://doi.org/10.1016/j.aim.2024.109524}.

\bibitem{Ball1986}
Keith Ball, \emph{Cube slicing in {$\mathbb R^n$}}, Proceedings of the American
  Mathematical Society \textbf{97} (1986), no.~3, 465--473.
\newblock \url{https://doi.org/10.1090/S0002-9939-1986-0840631-0}.

\bibitem{Ball1989}
\bysame, \emph{Volumes of sections of cubes and related problems}, Geometric
  Aspects of Functional Analysis (1987--88) (Joram Lindenstrauss and Vitali~D.
  Milman, eds.), Lecture Notes in Mathematics, vol. 1376, Springer, Berlin,
  1989, pp.~251--260.
\newblock \url{https://doi.org/10.1007/BFb0090058}.

\bibitem{Ball1991}
\bysame, \emph{Volume ratios and a reverse isoperimetric inequality},
Journal of the London Mathematical Society (2) \textbf{44} (1991), no.~2, 351--359.
\newblock \url{https://doi.org/10.1112/jlms/s2-44.2.351}.

\bibitem{Barthe1998}
Franck Barthe, \emph{On a reverse form of the {Brascamp--Lieb} inequality},
  Inventiones Mathematicae \textbf{134} (1998), no.~2, 335--361.
\newblock \url{https://doi.org/10.1007/s002220050267}.

\bibitem{BGMN2005}
Franck Barthe, Olivier Gu\'edon, Shahar Mendelson, and Assaf Naor, \emph{A
  probabilistic approach to the geometry of the {$\ell_p^n$}-ball}, The Annals
  of Probability \textbf{33} (2005), no.~2, 480--513.
\newblock \url{https://doi.org/10.1214/009117904000000874}.

\bibitem{BCCT2008}
Jonathan Bennett, Anthony Carbery, Michael Christ, and Terence Tao, \emph{The
  {Brascamp--Lieb} inequalities: finiteness, structure and extremals},
  Geometric and Functional Analysis \textbf{17} (2008), no.~5, 1343--1415.
\newblock \url{https://doi.org/10.1007/s00039-007-0619-6}.

\bibitem{BFH2019}
K.~J. B\"or\"oczky, F.~Fodor, and D.~Hug, \emph{Strengthened volume
inequalities for \(L_p\) zonoids of even isotropic measures},
Transactions of the American Mathematical Society \textbf{371}
(2019), no.~1, 505--548.
\newblock \url{https://doi.org/10.1090/tran/7299}.

\bibitem{BFH2025}
\bysame, \emph{Strengthened inequalities for the mean width and the $\ell$-norm of origin symmetric convex bodies},
Mathematische Annalen \textbf{393} (2025), 271--316.
\newblock \url{https://doi.org/10.1007/s00208-025-03228-0}.

\bibitem{BrascampLieb1976}
Herm~Jan Brascamp and Elliott~H. Lieb, \emph{Best constants in {Young}'s
  inequality, its converse, and its generalization to more than three
  functions}, Advances in Mathematics \textbf{20} (1976), no.~2, 151--173.
\newblock \url{https://doi.org/10.1016/0001-8708(76)90184-5}.

\bibitem{BE2026}
Silouanos Brazitikos and Alexandros Eskenazis,
\emph{The dimension-free threshold for maximal sections of $\ell_p$},
preprint (2026).
\newblock \url{https://arxiv.org/abs/2611.07619}.

\bibitem{CarlenCordero2009}
Eric~A. Carlen and Dario Cordero-Erausquin, \emph{Subadditivity of the entropy
  and its relation to {Brascamp--Lieb} type inequalities}, Geometric and
  Functional Analysis \textbf{19} (2009), no.~2, 373--405.
\newblock \url{https://doi.org/10.1007/s00039-009-0001-y}.

\bibitem{CNT2022}
Giorgos Chasapis, Piotr Nayar, and Tomasz Tkocz, \emph{Slicing {$\ell_p$}-balls
  reloaded: stability, planar sections in {$\ell_1$}}, The Annals of
  Probability \textbf{50} (2022), no.~6, 2344--2372.
\newblock \url{https://doi.org/10.1214/22-AOP1584}.

\bibitem{ENT2024Resilience}
Alexandros Eskenazis, Piotr Nayar, and Tomasz Tkocz,
\emph{Resilience of cube slicing in $\ell_p$},
Duke Mathematical Journal \textbf{173} (2024), no.~17, 3377--3412.
\newblock \url{https://doi.org/10.1215/00127094-2024-0004}.

\bibitem{Hensley1979}
Douglas Hensley, \emph{Slicing the cube in $\mathbb R^n$ and probability
  (bounds for the measure of a central cube slice in $\mathbb R^n$ by
  probability methods)}, Proceedings of the American Mathematical Society
  \textbf{73} (1979), no.~1, 95--100.
\newblock \url{https://doi.org/10.1090/S0002-9939-1979-0512066-8}.

\bibitem{IriyehShibata2020}
Hiroshi Iriyeh and Masataka Shibata, \emph{Symmetric Mahler's conjecture for the volume product in the 3-dimensional case},
Duke Mathematical Journal \textbf{169} (2020), no.~6, 1077--1134.
\newblock \url{https://doi.org/10.1215/00127094-2019-0072}.

\bibitem{Ivanov2021Cross}
Grigory Ivanov, \emph{On the volume of projections of the cross-polytope},
  Discrete Mathematics \textbf{344} (2021), no.~5, 112312.
\newblock \url{https://doi.org/10.1016/j.disc.2021.112312}.

\bibitem{Ivanov2021Frames}
\bysame, \emph{Tight frames and related geometric problems}, Canadian
  Mathematical Bulletin \textbf{64} (2021), no.~4, 942--963.
\newblock \url{https://doi.org/10.4153/S000843952000096X}.

\bibitem{IvanovTsiutsiurupa2021}
Grigory Ivanov and Igor Tsiutsiurupa, \emph{On the volume of sections of the
  cube}, Analysis and Geometry in Metric Spaces \textbf{9} (2021), no.~1,
  1--18.
\newblock \url{https://doi.org/10.1515/agms-2020-0103}.

\bibitem{KyENT2026}
Vu Khac Ky, \emph{A proof of the Eskenazis--Nayar--Tkocz conjecture on
entropy maximization for cube projections}, arXiv:2611.18798 (2026).
\newblock \url{https://arxiv.org/abs/2611.18798}.

\bibitem{Lieb1990}
Elliott~H. Lieb, \emph{Gaussian kernels have only {Gaussian} maximizers},
  Inventiones Mathematicae \textbf{102} (1990), 179--208.
\newblock \url{https://doi.org/10.1007/BF01233426}.

\bibitem{LivshytsPaourisPivovarov2016}
Galyna Livshyts, Grigoris Paouris, and Peter Pivovarov, \emph{On sharp bounds
  for marginal densities of product measures}, Israel Journal of Mathematics
  \textbf{216} (2016), no.~2, 877--889.
\newblock \url{https://doi.org/10.1007/s11856-016-1431-5}.

\bibitem{Lyons2003}
Russell Lyons, \emph{Determinantal probability measures}, Publications
  Math\'ematiques de l'IH\'ES \textbf{98} (2003), 167--212.
\newblock \url{https://doi.org/10.1007/s10240-003-0016-0}.

\bibitem{MMX2017}
Mokshay Madiman, James Melbourne, and Peng Xu, \emph{Rogozin's convolution inequality for locally compact groups},
preprint, arXiv:1705.00642 (2017).
\newblock \url{https://arxiv.org/abs/1705.00642}.

\bibitem{MelbourneRoberto2022}
James Melbourne and Cyril Roberto, \emph{Quantitative form of {Ball}'s cube
  slicing in {$\mathbb R^n$} and equality cases in the min-entropy power
  inequality}, Proceedings of the American Mathematical Society \textbf{150}
  (2022), no.~8, 3595--3611.
\newblock \url{https://doi.org/10.1090/proc/15944}.

\bibitem{MeyerPajor1988}
Mathieu Meyer and Alain Pajor, \emph{Sections of the unit ball of
  {$\ell_p^n$}}, Journal of Functional Analysis \textbf{80} (1988), no.~1,
  109--123.
\newblock \url{https://doi.org/10.1016/0022-1236(88)90068-7}.

\bibitem{NayarTkocz2023}
Piotr Nayar and Tomasz Tkocz, \emph{Extremal sections and projections of
  certain convex bodies: a survey}, Harmonic Analysis and Convexity (Alexander
  Koldobsky and Alexander Volberg, eds.), Advances in Analysis and Geometry,
  vol.~9, De Gruyter, Berlin, 2023, pp.~343--390.
\newblock \url{https://doi.org/10.1515/9783110775389-008}.

\bibitem{OpenAIMahler2026}
OpenAI, \emph{The symmetric Mahler conjecture and its equality cases},
preprint, September 22, 2026.
\newblock \url{https://github.com/openai/math/tree/main/preprints/The-symmetric-Mahler-conjecture-and-its-equality-cases-September-22-2026}.

\bibitem{Polya1913}
G. P\'olya, \emph{Berechnung eines bestimmten Integrals},
Mathematische Annalen \textbf{74} (1913), no.~2, 204--212.

\bibitem{Pooley2012}
Benjamin C. Pooley, \emph{Upper bounds on the cross-sectional volumes of
  cubes and other problems}, Master of Mathematics dissertation,
  University of Warwick, March 2012.
\newblock \url{https://warwick.ac.uk/fac/sci/maths/people/staff/bpooley/pooleymmathproject2012.pdf}.

\bibitem{Rogozin1987}
B.~A. Rogozin, \emph{Estimation of the maximum of a convolution of bounded
  densities}, Theory of Probability and its Applications \textbf{32} (1987),
  no.~1, 48--56.
\newblock \url{https://doi.org/10.1137/1132004}.

\bibitem{RV2015}
Mark Rudelson and Roman Vershynin, \emph{Small ball probabilities for linear images of high-dimensional distributions},
International Mathematics Research Notices \textbf{2015} (2015), no.~19, 9594--9617.
\newblock \url{https://doi.org/10.1093/imrn/rnu243}.

\bibitem{Sbai1987}
Sidi Driss Sbai Msahli, \emph{Sections vectorielles du cube et d'autres
boules unit\'e de $\mathbb R^n$}, doctoral thesis, Universit\'e des
sciences et techniques de Lille Flandres Artois, 1987.
\newblock \url{https://pepite-depot.univ-lille.fr/LIBRE/Th_Num/1987/55376-1987-23.pdf}.

\bibitem{Schneider}
Rolf Schneider, \emph{Convex bodies: The {Brunn--Minkowski} theory}, second
  expanded ed., Encyclopedia of Mathematics and its Applications, vol. 151,
  Cambridge University Press, Cambridge, 2014.
\newblock \url{https://doi.org/10.1017/CBO9781139003858}.

\bibitem{Talagrand1996}
Michel Talagrand, \emph{Transportation cost for {Gaussian} and other product
  measures}, Geometric and Functional Analysis \textbf{6} (1996), no.~3,
  587--600.
\newblock \url{https://doi.org/10.1007/BF02249265}.

\bibitem{Vaaler1979}
Jeffrey~D. Vaaler, \emph{A geometric inequality with applications to linear
  forms}, Pacific Journal of Mathematics \textbf{83} (1979), no.~2, 543--553.
\newblock \url{https://doi.org/10.2140/pjm.1979.83.543}.

\end{thebibliography}
\end{document}